Subject: A Noncommutative Chromatic Symmetric Function

Abstract:

In~\cite{stan}, Stanley associated with a graph $G$ a symmetric
function $X_G$ which reduces to $G$'s chromatic polynomial $\chr$
under a certain specialization of variables.  He then proved various
theorems generalizing results about $\chr$, as well as new ones that
cannot be interpreted on the level of the chromatic polynomial.
Unfortunately, $X_G$ does not satisfy a Deletion-Contraction Law which
makes it difficult to apply the useful technique of induction.  
We introduce a symmetric function $Y_G$ in noncommuting
variables which does have such a law and specializes to $X_G$ when the
variables are allowed to commute.  This permits us to further
generalize some of Stanley's theorems and prove them in a uniform
and straightforward manner.  Furthermore, we make some progress on the
({\bf 3+1})-free  Conjecture of Stanley and
Stembridge~\cite{stembridge}.

-----------------------------------------------------------------------
\documentclass[titlepage,12pt]{article}
\usepackage{amssymb}
\usepackage{latexsym}

\newtheorem{thm}{Theorem}[section]
\newtheorem{proposition}[thm]{Proposition}
\newtheorem{corollary}[thm]{Corollary}
\newtheorem{lemma}[thm]{Lemma}
\newtheorem{conjecture}[thm]{Conjecture}
\newtheorem{exa}[thm]{Example}
\newtheorem{question}[thm]{Question}

\def\fu{ X_{G} }
\def\su{\sum_{ \kappa } x_ {\kappa( v_{1} )} 
x_{\kappa(v_{2})} \cdots x_{ \kappa ( v_{d})}}
\def\ind{\hspace{-.05 in}\uparrow}
\def\st{\tilde{S}}
\def\sub{S \subseteq E}
\def\non{ Y_{G}}
\def\nond{  Y_{G - e}}
\def\nonc { Y_{G/e}}
\def\broke{B_{G}}
\def\n{\noindent}

\def\pd{$\Pi_{d}$ }
\def\p{\Pi_{d} }
\def\comp{\circ}
\def\chr{{\cal X}_{G}(n)}
\def\rep{Relabeling Proposition }
\def\del{Deletion-Contraction Proposition }
\def\v{v_{0}}

\def\plus{\pi+(d+1)}
\def\slash{\pi / d+1}
\def\mucmm{{m-1\choose \de_{1},\ldots,\de_{j}-1, 
\ldots,\de_{q+2}}}

\def\mucm{{m\choose \de_{1},\ldots,\de_{j}, 
\ldots,\de_{q+2}}}
\def\fact{{(-1)^{q}q!\over{\prod_{r=1}^{m+1}k_{r
}!}}}
\def\cont{\subseteq}
\def\pta{(\tau )}
\def\ptas{(\tau /d)}
\def\ptap{(\tau +(d))}
\def\ptapp{(\tau +(d)+(d+1))}

\def\ptaps{(\tau +(d)/d+1)}
\def\ptasp{(\tau/d,d+1)}
\def\r{\right)}
\def\l{\left(}
\def\h{\hspace{1mm}}

\newcommand{\ben}{\begin{enumerate}}
\newcommand{\een}{\end{enumerate}}
\newcommand{\ble}{\begin{lemma}}
\newcommand{\ele}{\end{lemma}}
\newcommand{\bth}{\begin{thm}}
\renewcommand{\eth}{\end{thm}}
\newcommand{\bpr}{\begin{proposition}}
\newcommand{\epr}{\end{proposition}}
\newcommand{\bco}{\begin{corollary}}
\newcommand{\eco}{\end{corollary}}
\newcommand{\bcon}{\begin{conj}}
\newcommand{\econ}{\end{conj}}
\newcommand{\bde}{\begin{defn}}
\newcommand{\ede}{\end{defn}}
\newcommand{\bex}{\begin{exa}}
\newcommand{\eex}{\end{exa}}
\newcommand{\barr}{\begin{array}}
\newcommand{\earr}{\end{array}}
\newcommand{\btab}{\begin{tabular}}
\newcommand{\etab}{\end{tabular}}
\newcommand{\beq}{\begin{equation}}
\newcommand{\eeq}{\end{equation}}
\newcommand{\bea}{\begin{eqnarray*}}
\newcommand{\eea}{\end{eqnarray*}}
\newcommand{\bce}{\begin{center}}
\newcommand{\ece}{\end{center}}
\newcommand{\bpi}{\begin{picture}}
\newcommand{\epi}{\end{picture}}
\newcommand{\bfi}{\begin{figure} \begin{center}}
\newcommand{\efi}{\end{center} \end{figure}}

\newcommand{\bsl}{\begin{slide}{}}
\newcommand{\esl}{\end{slide}}

\newcommand{\Qed}{\rule{1ex}{1ex} \medskip}

\newcommand{\ol}{\overline}

\newcommand{\sbe}{\subseteq}

\newcommand{\setm}{\setminus}

\newcommand{\iso}{\cong}
\newcommand{\con}{\equiv}

\newcommand{\mh}{\hat{0}}
\newcommand{\Mh}{\hat{1}}

\newcommand{\zh}{\hat{0}}

\newcommand{\ipr}[1]{\langle{#1}\rangle}

\newcommand{\ra}{\rightarrow}

\newcommand{\al}{\alpha}

\newcommand{\ga}{\gamma}
\newcommand{\de}{\delta}

\newcommand{\si}{\sigma}

\newcommand{\bbC}{{\mathbb C}}

\newcommand{\bbN}{{\mathbb N}}

\newcommand{\cS}{{\cal S}}

\newcommand{\St}{\tilde{S}}

\newcommand{\dil}{\displaystyle}

\newcommand{\ds}{\displaystyle}
\newcommand{\defn}{\stackrel{\rm def}{=}}

\newtheorem{definition}[thm]{Definition}

\def\fu{ X_{G} }
\def\su{\sum_{ \kappa } x_ {\kappa( v_{1} )} 
x_{\kappa(v_{2})} \cdots x_{ \kappa ( v_{d})}}
\def\ind{\hspace{-.05 in}\uparrow}
\def\st{\tilde{S}}
\def\sub{S \subseteq E}
\def\non{ Y_{G}}
\def\nond{  Y_{G \setminus e}}
\def\nonc { Y_{G/e}}
\def\broke{B_{G}}
\def\n{\noindent}

\def\pd{$\Pi_{d}$ }
\def\p{\Pi_{d} }
\def\comp{\circ}
\def\chr{{\cal X}_{G}(n)}
\def\rep{Relabeling Proposition }
\def\del{Deletion-Contraction Proposition }

\begin{document}
\pagestyle{empty}
\title{A Noncommutative Chromatic Symmetric 
Function}
\author{David D. Gebhard\\
Box 7281 \\
Lenoir-Rhyne College\\
Hickory, NC 28603\\[5pt]
and\\[5pt]
Bruce E. Sagan\\
Department of Mathematics\\ 
Michigan State University\\
East Lansing, MI 48824-1027\\
sagan@math.msu.edu}

\date{\today \\[1in]
	\begin{flushleft}
	\end{flushleft}
       }
\maketitle

\begin{flushleft} Proposed running head: 
\end{flushleft}
	\begin{center} 
Chromatic Symmetric Function
	\end{center}

Send proofs to:
\begin{center}
Bruce E. Sagan \\ Department of Mathematics 
\\Michigan State
University \\ East Lansing, MI 48824-1027\\[5pt]
Tel.: 517-355-8329\\
FAX: 517-432-1562\\
Email: sagan@math.msu.edu
\end{center}

	\begin{abstract}
In~\cite{stan}, Stanley associated with a graph $G$ a symmetric
function $X_G$ which reduces to $G$'s chromatic polynomial $\chr$
under a certain specialization of variables.  He then proved various
theorems generalizing results about $\chr$, as well as new ones that
cannot be interpreted on the level of the chromatic polynomial.
Unfortunately, $X_G$ does not satisfy a Deletion-Contraction Law which
makes it difficult to apply the useful technique of induction.  
We introduce a symmetric function $Y_G$ in noncommuting
variables which does have such a law and specializes to $X_G$ when the
variables are allowed to commute.  This permits us to further
generalize some of Stanley's theorems and prove them in a uniform
and straightforward manner.  Furthermore, we make some progress on the
({\bf 3+1})-free  Conjecture of Stanley and Stembridge~\cite{stembridge}.
	\end{abstract}
\pagestyle{plain}

\section{Introduction}

Let $G$ be a finite graph with verticies
$V=V(G)=\{v_1,v_2,\ldots,v_d\}$ and edge 
set $E=E(G)$.  We permit our graphs to have loops and multiple edges.
Let $\chr$ be the {\it chromatic polynomial} of $G$, i.e.,
the number of proper colorings $\kappa:V\ra\{1,2,\ldots,n\}$.  
({\it Proper} means that $vw\in E$ implies $\kappa(v)\neq\kappa(w)$.)

In~\cite{stan,stanley}, R.\ P.\ Stanley 
introduced a symmetric function, $ X_{G}$, 
which generalizes $\chr$ as follows.
Let $x=\{x_1,x_2,\ldots\}$ be a countably infinite set of commuting 
indeterminates.  Now define
    \[  X_{G} =  X_{G}(x_ {1},x_{ 2}, \ldots) = 
	\sum_{ \kappa } x_{ \kappa(v_{1})} \ldots x_{ \kappa (v_{d})}\]
where the sum ranges over all proper colorings,
    $\kappa: V(G) \rightarrow \{1,2,\ldots\}$.
It is clear from the definition that $\fu$ is a 
symmetric function, since permuting the colors of a proper coloring
leaves it proper, and is homogeneous of degree $d=|V|$.  Also the
specialization $X_{G}(1^{n})$  obtained by setting 
$x_{1}=x_{2}=\cdots =x_{n}=1$, and $x_{i}=0$  
for all $i>n$  yields ${\cal X}_{G}(n)$.

Stanley used $X_G$ to generalize various results about the chromatic
polynomial as well as proving new theorems that only apply to the
symmetric function.  However, there is a problem when trying to find
a deletion-contraction law for $X_G$.  To see what goes wrong, suppose 
that for $e\in E$ we let $G\setm e$ and $G/e$ denote $G$ with the
$e$ deleted and contracted, respectively.  Then $X_G$ and 
$X_{G\setm e}$ are homogeneous of degree $d$ while $X_{G/e}$ is
homogeneous of degree $d-1$ so there can be no linear relation
involving all three.  We should note that Noble and
Welsh~\cite{nw:wgp} have a deletion contraction method for computing
$X_G$ equivalent to~\cite[Theorem 2.5]{stan}.  However, it only works
in the larger category of vertex-weighted graphs and only for the
expansion of $X_G$ in terms of the power sum symmetric functions.
Since we are interested in other bases as well, we take a different
approach.  

In this paper we define an analogue, $Y_G$, of $\fu$ which is a symmetric
function in {\em  noncommuting} variables.  
(Note that  these noncommutative 
symmetric functions are different from the
noncommuting symmetric functions  studied by 
Gelfand and others, see \cite{thibon} for 
example.)  The reason for not letting the
variables commute is so that we can keep track of the color 
which $\kappa$ assigns to each vertex.
This permits us to prove a Deletion-Contraction Theorem for $Y_G$ and
use it to derive generalizations of results about $X_G$ in a
straightforward manner by induction, as well as make progress on a
conjecture. 

The rest of this paper is organized as follows.  In the next section
we begin with some basic background about
symmetric functions in noncommuting variables
(see also~\cite{me}).  In Section~\ref{YG} we define $Y_G$ and derive
some of its basic properties, including the Deletion-Contraction Law.
Connections with acyclic orientations are explored in Section~\ref{ao}.
The next three sections are devoted to making some progress on the
({\bf 3+1})-free  Conjecture of Stanley and Stembridge~\cite{stembridge}.
Finally we end with some comments and open questions.

\section{Noncommutative symmetric functions}
 
Our noncommutative symmetric functions will 
be indexed by elements of the partition lattice.  
We let \pd  \/ denote the lattice of set 
partitions $\pi$ of $\{1,2, \ldots,d\}:=[d]$, ordered 
by refinement. We write $\pi=B_1/B_2\ldots/B_k$ if 
$\uplus_i B_i=[d]$ and call $B_i$ a {\it block} of $\pi$.
The  meet (greatest lower bound) 
of the elements $\pi$ and $\sigma$  is denoted 
by $\pi \wedge \sigma$.  We use $\mh$ to  
denote the unique minimal element, and $\Mh$ for
the unique maximal element.

For   $\pi \in \p$ we define $\lambda(\pi)$ to 
be the integer partition of $d$ whose parts are 
the 
block sizes 
of $\pi$. 
Also, if $\lambda(\pi)= 
(1^{r_{1}},2^{r_{2}}, 
\ldots,d^{r_{d}})$,  we will need 
the constants  
\bea
|\pi|	&=&r_{1}!r_{2}! \cdots r_{d}! \mbox{ and }  \\
\pi!	&=&1!^{r_{1}}2!^{r_{2}} \cdots d!^{r_{d}}.
\eea

We now introduce the vector space for the 
noncommutative symmetric functions. Let 
$\{x_{1},x_{2},x_{3}, \ldots \}$ be a set of 
{\em noncommuting} variables.   We define the 
{\em 
noncommutative monomial symmetric functions}, 
$m_{\pi}$, by
\begin{equation}
 m_{\pi}= \sum_{i_{1},i_{2},\ldots ,i_{d}} 
x_{i_{1}}x_{i_{2}}\cdots x_{i_{d}}, 
\label{mono}
\end{equation} 
where the sum is over all sequences
$i_{1},i_{2},\ldots ,i_{d} $
of positive integers 
such that $i_{j}=i_{k}$ if and only if $j$ and 
$k$ are  in the same block of $\pi$.  For 
example, we get
$$m_{13/24}=x_{1}x_{2}x_{1}x_{2} 
+x_{2}x_{1}x_{2}x_{1}  +x_{1}x_{3}x_{1}x_{3}+ 
x_{3}x_{1}x_{3}x_{1} +\cdots  $$ 
for the partition $\pi = 13/24$.

{From} the definition it is easy to see that letting 
the $x_{i}$ commute transforms $m_{\pi}$ into 
$|\pi|m_{\lambda(\pi)}$, a multiple of the ordinary 
monomial symmetric function.  
The noncommutative monomial symmetric functions, 
$\{m_{\pi}:\pi \in \p, d \in {\bbN}\}$, are 
linearly independent 
over $\bbC$, and we  call 
their span 
the algebra of {\em noncommutative symmetric 
functions}.  
 
%
%
There are two other bases of this algebra that will interest us.  
One of them consists of the {\em 
noncommutative power sum symmetric functions} 
given by
\begin{equation}
 p_{\pi}\defn \sum_{\sigma \geq \pi}m_{\sigma} = 
 \sum_{i_{1},i_{2},\ldots ,i_{d}} 
x_{i_{1}}x_{i_{2}}\cdots 
x_{i_{d}},  
\label{power1}
\end{equation} 
where the second sum is over all positive integer sequences
$i_{1},i_{2},\ldots ,i_{d}$ 
such that $i_{j}=i_{k}$ 
if $j$ and $k$ are both in the same block of 
$\pi$.  The other basis contains the
{\em noncommutative 
elementary symmetric functions} defined by
\begin{equation}
e_{\pi}\defn\sum_{\sigma:\sigma \wedge \pi =\mh}m_{\sigma}= 
\sum_{i_{1},i_{2},\ldots,i_{d}} x_{i_{1}}x_{i_{2}}\cdots x_{i_{d}},
\label{elem1}
\end{equation}  
 where the second sum is over all sequences
$i_{1},i_{2},\ldots ,i_{d}$ of positive integers such that 
$i_{j}\neq i_{k}$ if $j$ and $k$ are both in the same block 
of $\pi$.  As an illustration of these definitions, we  see 
that 
\begin{eqnarray*}p_{13/24}&=&x_{1}x_{2}x_{1}x_{2
}+x_{2}x_{1}x_{2}x_{1} + \cdots 
+x_{1}^{4}+x_{2}^{4}+ \cdots\\
&=& m_{13/24}+m_{1234}\\
\end{eqnarray*}
and that
 \begin{eqnarray*}e_{13/24}&=&
x_{1}^{2}x_{2}^{2}+\cdots + x_{1}x_{2}^{2}x_{1} 
+\cdots +x_{1}^{2}x_{2}x_{3}+\cdots \\
& &  + x_{1}x_{2}^{2}x_{3}+\cdots + 
x_{1}x_{2}x_{3}^{2}+\cdots 
+x_{1}x_{2}x_{3}x_{1}+\cdots 
+x_{1}x_{2}x_{3}x_{4}\cdots\\
&=& 
m_{12/34}+m_{14/23}+m_{12/3/4}+m_{1/23/4}+m_{1/2
/34}+m_{14/2/3}+m_{1/2/3/4}.\\
\end{eqnarray*}
Allowing the variables to commute  transforms 
$p_{\pi}$ into $p_{\lambda(\pi)}$ and 
$e_{\pi}$ into 
$\pi!e_{\lambda(\pi)}$.
We  may also use these definitions to derive the  
change-of-basis formulae  found in the appendix of Doubilet's 
paper \cite{doubilet} which show
\begin{eqnarray}
m_{\pi}	&=&\dil\sum_{\sigma \geq \pi}\mu(\pi,\sigma)p_{\sigma},\\[10pt]
m_{\pi} &=&\dil\sum_{\tau \geq \pi}\frac{\mu(\pi,\tau)}{\mu(\mh,\tau)}
		\sum_{\sigma \leq \tau}\mu(\sigma,\tau)e_{\sigma},\\[10pt]
e_\pi	&=&\dil\sum_{\si\le\pi} \mu(\zh,\si) p_\si,\mbox{ and}
						\label{ep}\\[10pt] 
p_\pi	&=&\dil \frac{1}{\mu(\zh,\pi)} \sum_{\si\le\pi} \mu(\si,\pi)e_\si,
						\label{pe}
\end{eqnarray}
where $\mu(\pi,\si)$ is the M\"obius function of $\Pi_n$.

 It should be  clear that 
these noncommutative symmetric functions are 
symmetric in the 
usual sense, i.e., they are invariant under the 
usual symmetric group action on the variables.  
However, it 
will be useful to define a new action of the 
symmetric group on the noncommutative symmetric 
functions which 
permutes the {\em positions} of the variables.
%
%
%
For $\delta \in {\cal S}_{d}$, we define 
$$\delta \comp m_{\pi} \defn m_{\delta (\pi)},$$ 
where the action of 
$\delta \in 
S_{d}$ on a set partition of $[d]$ is the 
obvious one acting on the elements of the 
blocks.  It 
follows 
that for any $\delta$ this action induces a 
vector space isomorphism, since it merely 
produces a 
permutation 
of the basis elements.  Alternatively we can 
consider this action to be defined on the 
monomials so 
that 
\[\delta \comp (x_{i_{1}}x_{i_{2}} \cdots 
x_{i_{k}}) \defn 
x_{i_{\delta^{-1}(1)}}x_{i_{\delta^{-1}(2)}}
\cdots 
x_{i_{\delta^{-1}(k)}}\] and extend  linearly.

Utilizing the first characterization of this 
action, it follows straight from  definitions 
(\ref{power1}) 
and (\ref{elem1}) that $\delta \comp 
p_{\pi}=p_{\delta (\pi)}$ and $\delta \comp 
e_{\pi}=e_{\delta 
(\pi)}$.

\section{$Y_{G}$, The noncommutative version} \label{YG}

We begin by defining our main object of study, $Y_{G}$. 
\begin{definition} For any graph $G$ 
with vertices labeled  $ v_{1}, v_{2}, \ldots, 
v_{d}  $  {\em in a fixed order}, define 
     \[ Y_{G} = \su =\sum_\kappa x_\kappa , \] 
where again  the sum is over all proper colorings $\kappa$
of $G$, but the $ x_{i}$ are now {\em 
noncommuting} 
variables.
\end{definition}

As an example, for $P_{3}$, the path on three 
vertices with edge set $\{v_1v_2,v_2v_3\}$, we can calculate 
\begin{eqnarray*}Y_{P_{3}}&=& 
x_{1}x_{2}x_{1}+x_{2}x_{1}x_{2}+x_{1}x_{3}x_{1} 
+ \cdots + x_{1}x_{2}x_{3} 
+x_{1}x_{3}x_{2} +\cdots +x_{3}x_{2}x_{1} + 
\cdots \\
&=&m_{13/2}+m_{1/2/3}. \label{path}
\end{eqnarray*}
 
Note that if $G$ has  loops then 
this sum is empty and  
$\non=0$.  Furthermore, $\non$ depends not only 
on $G$, but also on the {\em labeling} of its 
vertices.   

In this section we will prove some results about the expansion of
$Y_G$ in various bases for the noncommutative symmetric functions and
show that it satisfies a Deletion-Contraction Recursion.  To obtain
the expansion in terms of monomial symmetric functions, note that
any partition $P$ of 
$V$ induces a set partition $\pi(P)$ of $[d]$ 
corresponding to the subscripts of the vertices. 
A partition $P$ of $V$ is {\it stable} if 
any two adjacent vertices are
in different blocks of $P$.  (If $G$ has a loop, 
there are no stable partitions.) The next result  follows directly
from the 
definitions.
\bpr					\label{exp}
We have
$$
\non =\sum_{P}m_{\pi(P)}
$$
where the sum is over all stable partitions, 
$P$, of $V$.\hfill\Qed
\epr

In order to show that $\non$ 
satisfies a Deletion-Contraction Recurrence 
it is necessary to have a distinguished edge.
Most of the time we will want this edge to be between the last two
vertices in the fixed order, but to permit an arbitrary edge choice
we will define an action of the 
symmetric group $\cS_d$
%
%
on a graph. For all 
$\delta 
\in {\cal S}_{d}$ we let $\delta$ act on the 
vertices of $G$ by 
$\delta(v_{i})=v_{\delta(i)}$.  This creates 
an 
action on graphs given by $\delta(G)=H$, where 
$H$ is just a relabeling of $G$.  
\bpr[Relabeling Proposition]  For any
graph $G$, we have
$$\delta \comp Y_{G} = Y_{\delta (G)},$$ where 
the vertex order $v_{1},v_{2},\ldots,v_{d}$ is 
used in both $Y_{G}$ and 
$Y_{\delta(G)}$. \label{relabel}
\epr
{\bf Proof.}
Let $\delta(G)=H$. We note that the action of 
$\delta$ produces a bijection between the stable 
partitions 
of $G$ 
and $H$.  Utilizing the previous proposition
and denoting the stable 
partitions of 
$G$ and $H$ by $P_{G}$ and $P_{H}$, respectively, 
we have 
\[ Y_{H}=\sum_{P_{H}}m_{\pi(P_{H})} = 
\sum_{P_{G}}m_{\delta(\pi(P_{G}))} = 
\sum_{P_{G}}\delta \comp 
m_{\pi(P_{G})}=\delta \comp 
\sum_{P_{G}}m_{\pi(P_{G})}=\delta \comp Y_{G}.  
\hspace{.3 in} \Qed \]

Using the \rep 
allows us, without loss 
of 
generality, to choose a labeling of $G$
%
%
with the distinguished edge for deletion-contraction being
$e=v_{d-1}v_{d}$.  
It is this edge for which we will derive the 
basic recurrence for 
$\non$.
\begin{definition} We define an operation 
called {\em induction},\hspace{.02 in} $\ind$, 
on the monomial  
\newline $x_{i_{1}}x_{i_{2}}  \cdots x_{i_{d-2}} 
x_{i_{d-1}}$,
by \[  (x_{i_{1}}x_{i_{2}}  \cdots x_{i_{d-2}} 
x_{i_{d-1}})\ind
      ~ = ~ x_{i_{1}}x_{i_{2}}  \cdots 
x_{i_{d-2}} x_{i_{d-1}}^{2}\]
and extend this operation linearly.
\end{definition}
 Note that 
this function takes a noncommutative symmetric 
function which is homogeneous of degree $d-1$ to 
one which is 
homogeneous of degree $d$. Context will make it 
clear whether  the word induction refers to this 
operation or to the 
proof technique.  
     
Sometimes we will also need to use induction on
an edge $e=v_{k}v_{l}$ so we extend the definition as follows. For
$k<l,$  define an operation~$\ind_{k}^{l}$  on noncommutative symmetric 
functions 
which simply repeats the variable in the $k^{\rm th}$ position again
in the $l^{\rm th}$.  That is, for a monomial 
     $x_{i_{1}}  \cdots x_{i_{k}} \cdots 
      x_{i_{d-1}}$,
define \[  (x_{i_{1}} \cdots x_{i_{k}} \cdots
          x_{i_{l-1}} x_{i_{l}} \cdots 
x_{i_{d-1}})\ind_{k}^{l}
       =x_{i_{1}} \cdots x_{i_{k}} \cdots 
x_{i_{l-1}}
          x_{i_{k}} x_{i_{l}} \cdots  
x_{i_{d-1}} \]
and extend linearly. 

Provided $G$ has an edge which 
is not 
a loop, 
we will usually  start by choosing a labeling such 
that $e=v_{d-1}v_{d}$. 
 We  also note here that if there is no such 
edge, then
\begin{equation}
\label{boundary}
Y_{G}= 
\left \{ \begin{array}{ll} p_{1/2/\cdots /d}=e_{1/2/\cdots /d} &
\mbox{ if }  G=\overline{K_{d}}\\ 0 & 
\mbox{ if } G \mbox{ has a loop,}\\
\end{array} \right.
\end{equation} where $\overline{K_{d}}$ is the 
completely disconnected graph on $d$ vertices.
We note that contracting an edge $e$ can create multiple 
edges (if there are vertices adjacent to both  of $e$'s endpoints) or
loops (if $e$ is part of a multiple edge), while
contracting  a loop deletes it.
%
%

\bpr[Deletion-Contraction Proposition]
For $e=v_{d-1}v_{d}$, we have
$$
\non = \nond \: -  \nonc \ind,
$$
where the  contraction of $e=v_{d-1}v_{d}$ is labeled 
$v_{d-1}$.
\epr

{\bf Proof.}  The proof is very similar to that 
of the Deletion-Contraction Property for ${\cal 
X}_{G}$.  We 
consider the proper colorings of $G- e$, which
 can 
be split disjointly into two types: \begin 
{enumerate}  
\item  proper colorings of $G- e$  with vertices 
$v_{d-1}$ and $v_{d}$ different colors;
\item proper colorings of $G- e$  with vertices 
$v_{d-1}$ and $v_{d}$ the same color.  \end 
{enumerate}

Those of the first type clearly correspond to 
proper colorings of $G$.   If $\kappa$ is a 
coloring of 
$G- e$ of the second type then (since the 
vertices $v_{d-1}$ and $v_{d}$ are the same 
color)  
we have    
      \[  x_ {\kappa( v_{1} )}  x_{ \kappa ( 
v_{2})} \cdots x_{ \kappa 
      (v_{d-1})}x_{\kappa(v_{d})} = 
      (x_ {\kappa( v_{1} )} 
x_{\kappa(v_{2})}\cdots x_{ \kappa ( v_{d-1})}) 
\ind =   x_ 
{\tilde{\kappa}} \ind 
\] where $\tilde {\kappa}$ is a proper coloring 
of $G/e$.  Thus we have  \( \nond = \non \: +  
\nonc \ind \).  
 \hfill $\Qed$
\vspace{.1in}

We note that if $e$ is a repeated edge, then 
the proper colorings of $G- e$ are exactly 
the same as those of $G$.  The fact that there 
are no proper colorings  of the second type 
corresponds to the 
fact that $G/e$ has at least one loop, and so it has no 
proper colorings.  Also note that because of our conventions for
contraction we always have
     \[ |E(G\setm e)|=|E(G/e)| = |E(G)|-1 \]
where $|\cdot|$ denotes cardinality.

It is easy to see how the operation of induction affects the monomial
and power sum symmetric functions.
For $\pi 
\in \Pi_{d-1}$ we  let  $\pi+(d) \in \p$ denote 
the partition $\pi$ with $d$ inserted into the 
block 
containing $d-1$.  From equations (\ref{mono}) 
and 
%
%
 (\ref{power1}) it is easy to see that 
$$m_{\pi}\ind = 
m_{\pi+(d)} \mbox{ ~ ~ and ~ ~ }  p_{\pi}\ind = 
p_{\pi+(d)}.$$  With this notation we can now 
provide an example of 
the Deletion-Contraction Proposition for 
$P_{3}$, where the vertices are labeled 
sequentially, and the 
distinguished edge is $e=v_{2}v_{3}$:
$$Y_{P_{3}}=Y_{P_{2}\uplus\{v_{3}\}} - 
Y_{P_{2}}\ind.$$  It is not difficult to compute 
 $$\begin{array}{ll} 
Y_{P_{2}}&= m_{1/2},\\
Y_{P_{2}}\ind & =m_{1/23},\\
Y_{P_{2}\uplus\{v_{3}\}} &=m_{1/2/3}+m_{1/23} 
+m_{13/2}.
\end{array}$$
This gives us  
\begin{eqnarray*}Y_{P_{3}}&=&m_{1/2/3}+m_{1/23}+
m_{13/2}-m_{1/23}\\
&=&m_{1/2/3}+m_{13/2},
\end{eqnarray*}
                                          which 
agrees with our previous calculation in equation (\ref{path}).
                                          
We may use the Deletion-Contraction Proposition to provide 
inductive proofs for noncommutative analogues   of some results of 
Stanley \cite{stan}.

\bth  \label{pexp}
For any  graph $G$,
 \[ \non = \sum_{\sub }(-1)^{|S|}p_{\pi (S)}, \] 
where $\pi(S)$ denotes the partition of $[d]$ 
associated with the vertex partition for  the 
connected components of the spanning subgraph of 
$G$ induced by $S$.
\eth

{\bf Proof.} We induct on the number of 
non-loops in $E$.  If $E$ consists only of 
$n$ loops, for $n \geq 0$, 
then for all $S\subseteq E(G)$, we will have 
$\pi(S)=1/2/\cdots /d$.  So
$$\sum_{\sub}(-1)^{|S|}p_{\pi (S)}= 
\sum_{\sub}(-1)^{|S|}p_{1/2/\ldots /d}=
\left\{ \begin{array}{ll} 
p_{1/2/\ldots /d} & \mbox{if $n=0$,}\\
   0 & \mbox{if $n>0.$}\\
  \end{array} \right.$$
This agrees with equation (\ref{boundary}).

Now, if $G$ has edges which are not loops, 
we use the \rep to obtain a labeling for $G$ 
with $e=v_{d-1}v_{d}$. From the \del   we know 
%
%
 that  
$\non = 
\nond \: -  \nonc \ind$ and by induction   
     \[\non= \sum_{S \subseteq E(G\setm e)}(-1)^{|S|}p_{\pi (S)} - 
        \sum _{\st \subseteq E(G/e)} 
(-1)^{|\st|} p_{\pi (\st)} \ind .\]

It should be clear that  
     \[ \sum_{S \subseteq E(G \setm e)} (-1)^{|S|} 
p_{\pi (S)} = 
        \sum_{  {\scriptstyle{ S \subseteq 
E(G)}} \atop{ \scriptstyle{e\not\in 
S}}}(-1)^{|S|}p_{\pi (S)}.\] 
Hence it suffices to show that 
\beq						\label{St}
 -  \sum _{\st \subseteq E(G/e)} 
(-1)^{|\st|} p_{\pi (\st)}\ind \hspace{.1 in} = 
        \sum_{ {\scriptstyle{S \subseteq E(G)}} 
\atop{\scriptstyle{e \in S}}}(-1)^{|S|}p_{\pi (S)} .
\eeq
To do so we define a map $ \Theta : 
\{\st\subseteq E(G/e) \}  \rightarrow  \{S 
\subseteq E(G) : e 
\in S \} $ by 
\[\Theta(\st) = S \mbox{, where $ S= \st \cup 
e$.}\]  
Then because of our conventions for contraction, $\Theta$ is a
bijection. 
Clearly $ \pi(S) = 
\pi(\tilde{S})+(d)$ giving $ 
p_{\pi(S)}=p_{\pi(\st)}\ind $. 
Furthermore $|S|=|\St|+1$ so
equation~(\ref{St}) follows and this
completes the proof. \hfill 
$ \Qed$ 
\vspace {.1in}

By letting the $x_{i}$ commute, we get
Stanley's Theorem 2.5 \cite{stan} as a 
corollary.  
Another results which we can obtain
by this method is 
%
%
 Stanley's generalization of   Whitney's Broken 
Circuit Theorem. 
 
  A  {\em circuit} is a closed walk,
$v_{1},v_{2},\ldots,v_{m},v_{1}$, with 
distinct vertices and edges.  Note that since we permit loops and
multiple edges, we can have $m=1$ or 2.  If we fix a 
total order on $E(G)$, a {\em broken circuit} is 
a circuit with its largest edge (with respect to 
the 
total order) removed.  Let $B_{G}$ denote the 
{\em broken 
circuit complex} of $G$, which is the set
 of all $S 
\subseteq E(G)$ which do {\em not} contain a 
broken circuit. 
Whitney's Broken Circuit Theorem states that the 
chromatic polynomial of a graph can be 
determined from its broken circuit complex.  
Before we prove our version of this theorem, 
however, we will need the 
following lemma, which appeared in the work of Blass and Sagan~\cite{bs}.  
\ble
For any non-loop $e$, there is a bijection 
between $\broke$ and $B_{G\setm e} \cup B_{G/e}$ 
given by 
\[ S \longrightarrow  \left \{ \begin{array}{ll} 
 \tilde{S}=S-{e} \in B_{G/e} & \mbox{if $e \in 
S$}\\
\tilde{S}=S \in B_{G\setm e} & \mbox{if $ e \notin S$,}\\
\end{array}
\right.
\]
where we take $e$ to be the first edge of $G$ in 
the total order on the edges . \hfill $\Qed$ \label{bijection}
\ele

Using this lemma, we can now obtain a 
characterization of $\non$ in terms of the 
broken circuit complex of $G$ for any fixed 
total 
ordering on the edges.

\bth 
We have \[ \non = \sum_{S\in \broke} 
(-1)^{|S|}p_{\pi (S)}, \] where $\pi(S)$ 
is as in Theorem~\ref{pexp}.
\eth

{\bf Proof.} We again induct on the number of 
non-loops in $E(G)$.  If the edge set consists 
only of $n$ loops, it 
should be clear that for $n>0$ we will have 
every edge being a circuit, and so the empty set 
is a broken circuit.  
Thus we have $$Y_{G}=\left\{ \begin{array}{ll} 
\sum_{S \in \{\phi\}} (-1)^{|S|}p_{\pi(S)}
=p_{1/2/\ldots /d} 
& \mbox{ if $n=0$,}\\
\sum_{S \in \phi}(-1)^{|S|}p_{\pi(S)} = 0 & 
\mbox{ if $n>0$,}
                                 \end{array} 
\right.$$
which matches equation (\ref{boundary}).

For $n>0$ and $e$ a non-loop, we  consider  
$\non = \nond \: -  \nonc \ind,$ and again apply 
induction.   From Lemma \ref{bijection} and 
arguments 
as in Proposition \ref{pexp}, we have
\[ \sum_{{\scriptstyle{S \in B_{G}}} \atop 
{\scriptstyle{e\notin S}}} (-1)^{S}p_{\pi (S)}= 
\sum_{{S\in B_{G\setm e}}} 
(-1)^{S}p_{\pi (S)} \] and 
%
%
\[ \sum_{ {\scriptstyle{S \subseteq E(G)}} 
\atop{\scriptstyle{e \in S}}}(-1)^{|S|}p_{\pi 
(S)} = -\sum _{\st \in B_{G/e}} 
(-1)^{|\st|} 
p_{\pi (\st)}\ind, \] which gives the result. 
\hfill $ \Qed$

\section{Acyclic orientations}			\label{ao}

An {\it orientation} of $G$ is a digraph obtained by assigning a
direction to each of its edges.  The orientation is {\it acyclic} if it
contains no circuits. 
A {\it sink} of an orientation is a vertex $v_0$ such that every
edge of $G$ containing it is oriented towards $v_0$.
There are some interesting results which 
relate the 
chromatic polynomial of a graph to the number of 
acyclic 
orientations of the graph and  the sinks of 
these  orientations.     The one which is the main motivation for this
section is the following theorem of Greene and
Zaszlavsky~\cite{gandz}.  To state it, we adopt the convention that 
the coefficient of $n^i$ in $\chr$ is $a_i$.
  
\bth \label{help} For any fixed 
vertex $v_{0}$, the 
number of acyclic orientations of $G$ with a 
unique 
sink at 
$v_{0}$ is $|a_{1}|$.\hfill\Qed
\eth

The original proof of this theorem uses the theory of hyperplane
arrangements.  For elementary bijective proofs, see~\cite{orient}.
Stanley~\cite{stan} has a stronger version of this result.

\bth \label{sta}
If $X_{G}=\sum_{\lambda}c_{\lambda}e_{\lambda}$, 
then the number of acyclic orientations of $G$ 
with $j$ sinks is 
given by $\ds{\sum_{l(\lambda)=j}c_{\lambda}.} 
\hfill \Qed$
\eth

We can prove an analogue of this theorem in the 
noncommutative setting by using techniques 
similar to his, but  
have not been able to do so using induction.  We 
can, however, inductively demonstrate a  related result
which, unlike Theorem~\ref{sta} implies Theorem~\ref{help}.  For this
result we need a lemma 
from~\cite{orient}.  To state it,
we denote the set of acyclic orientations of $G$ 
by  ${\cal A}(G)$, and the set of acyclic 
orientations of $G$ with 
a unique sink at $v_{0}$ by ${\cal A}(G,v_{0})$. 
For completeness, we provide a proof.

\ble \label{useful1}For any fixed vertex 
$v_{0}$, and any edge $e=uv_{0}, u\neq v_{0}$, 
the map 
$$D \longrightarrow \left \{  \begin{array}{ll}
D\setm e\in {\cal A}(G\setm e,v_{0}) & \mbox{ if $D\setm e 
\in{\cal A}(G\setm e,v_{0})$}\\
D/e\in {\cal A}(G/e,v_{0})  & \mbox{ if $D\setm e 
\notin{\cal A}(G\setm e,v_{0}),$} 
\end{array} \right. $$
 is a bijection between ${\cal A}(G,v_{0})$ and 
${\cal A}(G\setm e,v_{0}) \uplus {\cal 
A}(G/e,v_{0})$, where  the vertex 
of $G/e$ formed by contracting $e$ is labeled 
$v_{0}$. \ele
{\bf Proof.}  
We must first show that this map is
well-defined, i.e., that in both  cases we obtain an acyclic
orientation with unique sink at $v_{0}$.  This
is true in the first case by definition.  In
case two, where $D\setm e
\notin{\cal A}(G\setm e,v_{0})$, it must be true that
$D\setm e$ has sinks both at $u$ and at $v_{0}$
(since deleting a
directed edge of $D$ will neither disturb the acyclicity
of the orientation nor cause the sink
at $v_{0}$ to be lost).  
Since $u$ and $v_{0}$
are the only sinks
in $D\setm uv_{0}$ the contraction  must have a
unique sink at $v_{0}$, and there will be no 
cycles formed.   Thus the orientation $D/e$ will be
in ${\cal A}(G/e,v_{0})$ and
this map is well-defined.

To see that this is a bijection, we exhibit the inverse.  This  is
obtained by simply orienting
all edges of $G$ as in $D\setm uv_{0}$ or $D/uv_{0}$
as appropriate, and then adding in the oriented edge
$\overrightarrow{uv_{0}}$.  Clearly this  map is also well-defined. 
\hfill  $\Qed$

We can now apply deletion-contraction to obtain a noncommutative
version of Theorem~\ref{help}.
\bth \label{help2}
Let $\non =\displaystyle{ \sum_{{ \pi \in 
\Pi_{d}}}c_{\pi}e_{\pi}}$. Then for any fixed 
vertex, 
$v_{0}$, 
  the number of acyclic orientations 
of $G$ with a unique sink at $v_{0}$ is 
$(d-1)!c_{[d]}$.
\eth
 {\bf Proof.}  We again induct on the number of 
non-loops in $G$.  In the base case, if all the 
edges of $G$ are 
loops, then $$Y_{G}= \left \{ \begin{array}{ll} 
e_{1/2/\ldots/d} &  \mbox{ if $G$ has no 
edges}\\
                                    0 & \mbox{ 
if $G$ has loops.}\\
                                    \end{array}
                                    \right.$$
 So $$c_{[d]}=\left \{ \begin{array}{ll} 1 & 
\mbox{ if $G=K_{1}$}\\
                                        0 & 
\mbox{ if $d>1$ or $G$ has loops}\\
                                        
\end{array}
                                        \right 
\} = |{\cal A}(G,\v)|.$$                                    

 If $G$ has non-loops, then by the \rep  we 
may let 
$e=v_{d-1}v_{d}$.  We 
know that $\non = \nond - \nonc \ind$.  Since we 
will 
only be interested in the leading coefficient, 
let 
\[ \non=ae_{[d]} + \sum_{{\sigma < 
[d]}}a_{\sigma}e_{\sigma},\]
\[\nond=be_{[d]} + \sum_{{\sigma < 
[d]}}b_{\sigma}e_{\sigma},\] and
\[\nonc=ce_{[d-1]} + \sum_{{\sigma < [d-1]}} 
c_{\sigma}e_{\sigma}\]
where $\leq$ is the partial order on set 
partitions.
By induction and Lemma \ref{useful1}, it is 
enough to show that 
$(d-1)!a=(d-1)!b+(d-2)!c.$

Utilizing the change-of-basis formulae~(\ref{ep}) and~(\ref{pe}) 
as well as the fact that for $\pi \in \Pi_{d-1}$ we have $p_{\pi}\ind = 
p_{\pi + (d)}$, we obtain
\begin{equation}
 e_{\pi} \ind = \sum_{{ \sigma \leq \pi}} {{\mu 
(\mh, 
\sigma)} \over{ \mu (\mh, \sigma 
+(d))}}\sum_{{\tau 
\leq 
\sigma +(d)}} \mu(\tau, \sigma +(d)) e_{\tau}.
\label{changeup}
\end{equation}  
With this formula, we compute the coefficient of 
$e_{[d]}$ from $\nonc \ind$. The only term which 
contributes comes from $ce_{[d-1]} \ind$, which 
gives us 

 \begin{eqnarray*}
ce_{[d-1]} \ind&=& c\sum_{\sigma \in \Pi_{d-1}}  
{{\mu 
(\mh, \sigma)} \over{ \mu (\mh, \sigma 
+(d))}}\sum_{{\tau 
\leq \sigma +(d)}} \mu(\tau, \sigma +(d)) 
e_{\tau}\\
 &=& c{{\mu(\mh, [d-1])}\over{\mu (\mh, 
[d])}}e_{[d]} + 
\sum_{\tau <[d]} d_{\tau}e_{\tau}\\
 &=&{-c \over{d-1}}e_{[d]}+\sum_{\tau <[d]} 
d_{\tau}e_{\tau}\\
\end{eqnarray*} 
Thus, from  $\non = \nond - \nonc \ind$ we have 
that 
\begin{eqnarray*} (d-1)!a&=&(d-1)!b+(d-1)!{c 
\over{d-1}}\\
 &=&(d-1)!b+(d-2)!c,
\end{eqnarray*} 
which completes the proof. \hfill  $\Qed$

This result implies Theorem~\cite{gandz} since under the specialization
$x_{1}=x_{2}=\cdots =x_{n}=1$, and $x_{i}=0$  
for $i>n$, $e_\pi$ becomes
$$
\prod_{i=1}^k n(n-1)(n-2)\cdots(n-|B_i|+1)
$$
where $\pi=B_1/B_2/\ldots/B_k$.  So if $k\ge2$ this polynomial in $n$
is divisible by $n^2$. Thus the only summand contributing to the linear
term of $\chi_G(n)$ is when $\pi=[d]$ and in that case the
coefficient has absolute value $(d-1)!c_{[d]}$.

The next corollary follows easily from the previous result.

\bco
If $\non=\displaystyle{ \sum_{{ \pi \in 
\Pi_{d}}}c_{\pi}e_{\pi}}$, then the number of 
acyclic 
orientations of 
$G$ with one sink is $d!c_{[d]}$.  \hfill $\Qed$
\eco

 
\section{Inducing $e_{\pi}$}

We now turn our attention to the expansion of 
$Y_{G}$ in terms of the 
elementary symmetric function basis.  We recall  
that for 
any fixed $\pi \in \Pi_{d}$ we use $\plus$  to 
denote the partition of $[d+1] $ formed
by inserting the element $(d+1)$ into the block 
of $\pi$ which contains 
$d$.  We will denote the block of $\pi$ which 
contains $d$ by $B_{\pi}$.
  We also let $\slash$ be the partition of 
$[d+1]$ formed by adding the 
block $\{d+1\}$ to $\pi$.  

It is necessary for us to understand 
the coefficients arising in 
$e_{\pi}\ind$ if we want to understand the 
coefficients of $Y_{G}$ which occur in its 
expansion in terms of the elementary
 symmetric function basis.  We have seen in  equation 
(\ref{changeup}) that the expression for 
$e_{\pi}\ind$ is rather complicated.  
 However, if the terms in the expression of 
$e_{\pi}\ind$ are grouped properly,  the 
coefficients in many of the groups 
will sum to zero.   Specifically,  we need to 
combine the coefficients from set 
partitions which are of the same type 
(as integer partitions), and whose block 
containing $d+1$ have the same size.  
Keeping track of the size of the block containing
$d+1$ will allow us to use deletion-contraction repeatedly.
%
%
%
%
To do this formally, we introduce a bit of notation.
Suppose $\al=(\al_1,\al_2,\ldots,\al_l)$ is a {\it composition}, i.e.,
an ordered integer partition.
Let $P(\al)$ be the set of all
partitions $\tau=B_1/B_2/\ldots/B_l$ of $[d+1]$  such that
\ben
\item $\tau\le\plus$, 
\item $|B_i|=\al_i$ for $1\le i\le l$, and
\item $d+1\in B_1$.
\een
The proper grouping for the 
terms of $e_{\pi}\ind$  is given by the 
following lemma.

\ble                                             
        \label{main}
If $\ds{e_{\pi}\uparrow ~= \sum_{\tau \in 
\Pi_{d+1}}c_{\tau}e_{\tau}}$, 
then $c_{\tau}=0$ unless $\tau \leq \plus$, and 
for any composition $\al$, we have 
$$ \sum_{\tau \in P(\al)}c_{\tau}=
\left \{ \begin{array}{ll}
1/|B_{\pi}| & \mbox{ if  $P(\al)=\{\slash\}$}, 
\\[2pt]
-1/|B_{\pi}| & \mbox{ if  $P(\al)=\{\plus\}$}, 
\\[2pt]
0 & \mbox{ else}.
\end{array}
\right.
$$
\ele

{\bf Proof.}  Fix $\pi \in \Pi_{d}$.  By 
equation (\ref{changeup})   
$$ e_{\pi}\ind = \sum_{\si \leq \pi} 
{\mu(\hat{0},\si)\over
{\mu(\hat{0}, \si +(d+1))}}\sum_{\tau \leq \si 
+(d+1)}\mu(\tau, \si +(d+1))e_{\tau}.$$
Hence we may express 
$$e_{\pi}\ind = \sum_{\tau \leq 
\plus}c_{\tau}e_{\tau},$$
where for any fixed $\tau \leq \plus$ we have
\begin{equation}          
c_{\tau}=\sum_{{\scriptstyle{\si \leq \pi} 
\atop{ \scriptstyle{\si +(d+1) \geq \tau} }}}
{-1 \over{|B_{\si}|}}\mu\l \tau, \si +(d+1)\r .
\label{ctau}
\end{equation}

We first note that if $\tau  = \slash \in 
P(\alpha)$, then $|P(\alpha)|=1$ and we have
 the interval $[\tau, \plus]\iso \Pi_{2}$. A 
simple computation shows that $c_{\slash}= 
1/|B_{\pi}|$.  Similarly,
 if $\tau  = \plus \in P(\alpha)$, then again  
$|P(\alpha)|=1$ and we can easily compute 
$c_{\plus}= -1/|B_{\pi}|$.

We now fix $\tau = B_{1}/B_{2}/\cdots / B_{q+2} 
/ \cdots /B_{l} \in 
P(\al)$ and without loss of generality we can
 let $ B_{1},B_{2},\cdots , B_{q+2} $ where $q\ge-1$ be the 
blocks of $\tau$ 
which are contained in $B_{\pi+(d+1)}$.  For 
notational convenience, we will  
also let $|B_{\pi+(d+1)}|=m+1$, where $m \geq 
1$.
Finally, let $\beta$ denote the partition 
obtained from $\tau$ by 
merging the blocks of $\tau$ which contain $d$ 
and $d+1$, allowing $\beta = \tau$ if $d$ and 
$d+1$ are in the same 
block of 
$\tau$. Replacing $\si+(d+1)$ by $\si \in 
\Pi_{d+1}$ in equation (\ref{ctau}), we  see 
that 
$$c_{\tau} = 
\sum_{\beta \leq \si \leq \plus} 
{-1\over{|B_{\si}|-1}}\mu(\tau, \si).$$  

Now for any $B \cont [d+1]$ we will consider the 
sets 
$$L(B)=\{\si \in 
\Pi_{d+1}: \{d,d+1\} \cont B \in \sigma , \mbox{ 
where }\beta \leq \si \leq \plus\}.$$  
The nonempty $L(B)$ partition the interval 
$[\beta , 
\plus]$ according to the content of the block 
containing $\{d, d+1\}$
 and so we may express 
$$c_{\tau}=\sum_{B}{-1\over{|B|-1}}\sum_{\si \in 
L(B)} 
\mu(\tau, \si).$$
To compute the inner sum, we need to consider 
the following 2  cases.

{\bf Case 1)} For some $k>q+2,~ B_{k}$ is 
strictly contained in a block
 of $\plus$.  
In this case, we see that each non-empty $L(B)$ 
forms a non-trivial
cross-section of a product of 
partition lattices, and so for this case
 $$\sum_{\si \in L(B)} \mu(\tau, \si)=0.$$
Thus these $\tau$ will not contribute 
to $\ds{\sum_{\tau \in P(\alpha)}c_{\tau}}.$

\vspace{.2 in}

{\bf Case 2)} For all $k>q+2, ~B_{k}$ is a block 
of $\plus$.  So, by abuse of notation, we can write
$\tau=B_1/\ldots/B_{q+2}$ and $P(\al)=P(\al_1,\ldots,\al_{q+2})$.
Also in this case, we can assume $q \geq 0$, since 
 we have already computed this sum when
$\tau = \plus$. 
%
%
Then we will show 
\begin{equation} \label{cases} 
 \frac{1}{|B|-1}\sum_{\si \in L(B)} \mu(\tau, 
\si)=
\left \{ \begin{array}{ll}
\frac{(-1)^{q+1}(q+1)!}{m} & \mbox{ if 
$B=B_{\plus}$,}\\[2pt]
\frac{(-1)^{q}q!}{m-\alpha_{i}} & \mbox{ if 
$B=B_{\plus}\setminus B_{i}$, \hspace{.01in} $2 
\leq i \leq q+2$,}\\[2pt]
0 & \mbox{ else.}
 \end{array}
\right.
\end{equation}
Indeed, it is easy to see that if $B=B_{\plus}$ then 
$L(B)=\{\plus\}$ and so this part is clear.
Also, if $B=B_{\plus}\setminus B_{i}$ for some 
$2\leq i \leq q+2$, then we
have $|L(B)|=1$ again and $\sum_{\si \in 
L(B)}
 \mu(\tau, \si)=(-1)^{q}q!.$
Otherwise, $L(B)$ again forms a non-trivial
cross-section of a product of partition 
lattices, and again gives us no 
net 
contribution to the sum.

We notice that since $\{d,d+1 \} \subseteq B$, the 
second case in (\ref{cases}) will only occur if 
$d\in B_{j}$ for 
$j\neq i$.  Adding up all these contributing terms gives us 
$$c_{\tau}=(-1)^{q}q!\left(\sum_{{\scriptstyle{i
=2} \atop{\scriptstyle{i \neq j}}}}^{q+2} 
{1\over{m-\al_{i}}} - {q+1\over{m}}\right).$$

In order to compute the sum over all $\tau \in 
P(\al)$, it will be convenient to consider all 
possible orderings for 
the block of $\tau$ containing $d$.
So for  $1\leq j \leq q+2$, let 
$$P(\al,j)=
\{ (B_{1},B_{2}, \ldots , B_{q+2})
\hspace{.1 in}| \hspace{.1 in}  
B_{1}/B_{2}/ \ldots /B_{q+2} \in P(\alpha),\ d\in B_{j}\}.$$  
The sequence $(B_{1},B_{2}, \ldots , B_{q+2})$ forms the {\em ordered}
set partition $\tau$.
%
%
We also define
$$ \de_{j} = \left\{ \begin{array}{ll}
\al_{j}-1 & \mbox{if $j=1$}\\
\al_{j} & \mbox{else},\\
\end{array}
\right. $$ 
so 
$$ |P(\al,j)|=\mucmm .$$  Thus we can see that 
$$\sum_{\tau \in P(\al,j)}c_{\tau}=\mucmm 
(-1)^{q}q!\left(\sum_{{\scriptstyle{i=2} 
\atop{\scriptstyle{i \neq j}}}}^{q+2} 
{1\over{m-\al_{i}}} - {q+1\over{m}}\right).$$
To obtain the sum over all $\tau \in P(\al)$ we
 need to sum over all $ P(\al,j)$ for  $1\leq j 
\leq q+2$.    However, if we let
 $k_{r}$  be the number of blocks $B_{i}, 1\leq 
i \leq q+2$, which
 have size $r$, then in the sum over all $ 
P(\al,j)$,  each
 $\tau \in P(\alpha)$ appears 
$\Pi_{r=1}^{m+1}k_{r}!$ times.
Combining all this information, we see that 
$$\sum_{\tau \in 
P(\al)}c_{\tau}=\fact\sum_{j=1}^{q+2}\mucmm\left
(\sum_{{\scriptstyle{i=2} 
\atop{\scriptstyle{i \neq j}}}}^{q+2} 
{1\over{m-\de_{i}}} - {q+1\over{m}}\right).$$

Hence it suffices to show that 
$$\sum_{j=1}^{q+2}\mucmm\left(
\sum_{{\scriptstyle{i=2} 
\atop{\scriptstyle{i \neq j}}}}^{q+2} 
{1\over{m-\de_{i}}} - {q+1\over{m}}\right)=0.$$
Using the multinomial recurrence we have,
$$\sum_{j=1}^{q+2}
\mucmm = \mucm$$
%
%
%
\n and so we need only show that 
$$ \sum_{j=1}^{q+2}\mucmm 
\sum_{{\scriptstyle{i=2} \atop{\scriptstyle{i 
\neq j}}}}^{q+2} 
{1\over{m-\de_{i}}}={q+1\over{m}}\mucm.$$

However, we may express  \begin{eqnarray*}
\sum_{j=1}^{q+2}\mucmm \sum_{{\scriptstyle{i=2} 
\atop{\scriptstyle{i \neq 
j}}}}^{q+2}{1\over{m-\de_{i}}} &=& 
\sum_{j=1}^{q+2}\frac{\mucm\de_{j}}{m}
\sum_{{\scriptstyle{i=2} \atop{\scriptstyle{i \neq 
j}}}}^{q+2} {1\over{m-\de_{i}}} \\  
&=& 
\frac{\mucm}{m}\sum_{j=1}^{q+2}
\sum_{{\scriptstyle{i=2} 
\atop{\scriptstyle{i \neq 
j}}}}^{q+2}{\de_{j}\over{m-\de_{i}}}
\\
&=&\frac{\mucm}{m}\sum_{i=2}^{q+2}{1\over{m-\de_
{i}}}
\sum_{{\scriptstyle{j=1} 
\atop{\scriptstyle{j\neq i}}}}^{q+2}\de_{j} 
\\
&=&\frac{\mucm}{m}\sum_{i=2}^{q+2}{1\over{m-\de_
{i}}}(m-\de_{i})\\
&=&{q+1\over{m}}\mucm.  \end{eqnarray*} 

\vspace{-.2 in}\hfill $\Qed$

\section{Some $e$-positivity results}

We wish to use Lemma~\ref{main} to prove some 
positivity theorems about 
$Y_{G}$'s  expansion in the elementary symmetric 
function basis.  If the coefficients of the 
elementary symmetric 
functions 
in this expansion are all non-negative, then we 
say that $Y_{G}$ is {\em $e$-positive}.
  Unfortunately, even for some of the simplest 
graphs, $Y_{G}$ is usually not 
 $e-$positive.  The only graphs which are 
obviously  $e-$positive
 are  the complete graphs on $n$ vertices and 
their complements,  for 
which we have 
$Y_{K_{n}}=e_{[n]}$ and  $Y_{{\overline 
K_{n}}}=e_{1/2/\cdots /n}$.  Even  paths, with 
the vertices labeled sequentially, are not  
$e-$positive, for we 
can compute that 
$Y_{P_{3}}={1\over{2}}e_{12/3}{-{1\over{2}}}e_{13/
2}+{1\over{2}}e_{1/23}
 +{1\over{2}}e_{123}$.
However, in this example we can see that 
while $Y_{P_{3}}$ is not 
$e-$positive, if we identify all the  
terms having the same type and the same size block 
 containing 3, the sum  will be non-negative for each
 of these sets. 


  This observation along with the proof of the 
previous lemma inspires us to define equivalence 
classes reflecting 
the 
sets $P(\alpha)$.  If the block of 
$\sigma$ containing $i$ is $B_{\sigma,i}$ and 
the block of $\tau$ 
containing 
$i$ is $B_{\tau,i}$, we define $$\sigma \con_{i} 
\tau \mbox{ iff $\lambda(\sigma)=\lambda(\tau)$ 
and 
$|B_{\sigma,i}|=|B_{\tau,i}|$}$$ 
and extend this definition so that $$e_{\sigma} 
\con_{i} e_{\tau} \mbox{ iff $\sigma \con_{i} 
\tau$}.$$
We let $(\tau)$ and $e_{(\tau)}$ denote the 
equivalence classes of $\tau$  and $e_{\tau}$,
respectively. Taking 
formal 
sums of these equivalence classes allows us to write 
expressions such as 
$$\sum_{\sigma \in \Pi_{d}}c_{\sigma}e_{\sigma} 
\con_{i} \sum_{\pta \subseteq \Pi_{d} 
}c_{\pta}e_{\pta} \mbox{ where 
$\ds{c_{\pta}=\sum_{\si \in \pta}c_{\si}}$}.$$  
We will refer to this equivalence relation as 
{\em congruence modulo 
$i$}.

Using this notation, we have $Y_{P_{3}}\con_{3} 
{1\over{2}}e_{(12/3)} 
+{1\over{2}}e_{(123)}$, since 
$e_{13/2}\con_{3}e_{1/23}$. We will say that a 
labeled graph $G$ (and similarly 
$Y_{G}$) is 
$(e)-${\em positive} 
if all the $c_{\pta}$ are non-negative 
for some labeling of $G$ 
and  suitably chosen 
congruence. 
We notice that the expansion of $Y_{G}$ for a 
labeled graph  may have all non-negative 
amalgamated coefficients 
for congruence modulo $i$, but not for congruence 
modulo $j$.  However, if a different labeling
for an $(e)$-positive graph is chosen, then we 
can always find a corresponding
congruence class to again  see
$(e)$-positivity.  This should be clear from the 
Relabeling 
Proposition.  

We now 
turn our attention to showing that  paths, 
cycles, and complete graphs with one edge 
deleted are all
 $(e)$-positive.    We begin with a few more 
preliminary results about this congruence 
relation and how it affects 
our 
induction of $e_{\pi}$.  
 
We note  that in the proof of Lemma \ref{main}, 
the roles played by the elements $d$ and $d+1$ 
are essentially 
interchangeable. That is, if we let 
${\tilde{P}(\alpha)}$ be the set of all 
partitions $\tau=B_1/B_2/\ldots/B_l$ of $[d+1]$  such that
\ben
\item $\tau\le\plus$, 
\item $|B_i|=\al_i$ for $1\le i\le l$, and
\item $d\in B_1$,
\een
and let $\tilde{\pi}$
be the partition $\pi\in \Pi_d$ with $d$ replaced by 
$d+1$, then  the same 
proof 
will show that 
$$ \sum_{\tau \in {\tilde{P}(\al)}}c_{\tau}=
\left \{ \begin{array}{ll}
1/|B_{\pi}| & \mbox{ if  
$\tilde{P}(\al)=\{\tilde{\pi}/d\}$}, \\[2pt]
-1/|B_{\pi}| & \mbox{ if  
$\tilde{P}(\al)=\{\tilde{\pi}+(d)\}$}, \\[2pt]
0 & \mbox{ otherwise.}
\end{array}
\right.
$$
Note that here $\tilde{\pi}+(d)$ is the 
partition obtained from $\tilde{\pi}$ by 
inserting the element $d$ into the 
block of 
$\tilde{\pi}$ containing $d+1.$
 This allows us to state a corollary in terms of 
the congruence relationship just defined.
 

\bco \label{dord+1}
 If $b=|B_{\pi}|$, then for any $\pi \in 
\Pi_{d}$, we have
 $$e_{\pi}\ind 
\con_{d+1}\frac{1}{b}e_{(\pi/d+1)} - 
\frac{1}{b}e_{(\pi+(d+1))} $$
and 
 $$e_{\pi}\ind 
\con_{d}\frac{1}{b}e_{(\tilde{\pi}/d)} - 
\frac{1}{b}e_{(\tilde{\pi}+(d))}.$$   
\vspace{-.2 in} \hfill $\Qed$
\eco

The next lemma simply verifies that the 
induction operation respects the congruence 
relation and follows immediately from equation~(\ref{changeup})
or the previous corollary.

\ble
If $e_{\gamma} \con_{d} e_{\tau}$, then $e_{\gamma}\ind 
\con_{d+1} e_{\tau}\ind. \hfill \Qed$
\ele

{From} this  we can extend induction to congruence 
classes in a well-defined manner: 
$$\mbox{ if ~ }
e_{\pi}\ind = \sum_{\tau \in 
\Pi_{d+1}}c_{\tau}e_{\tau}  \mbox{ ~ then  ~ } 
e_{(\pi)}\ind \con_{d+1} \sum_{(\tau) \subseteq 
\Pi_{d+1}}c_{(\tau)}e_{(\tau)}.
$$

In order to use induction to prove the  
$(e)$-positivity of a graph $G$, we will usually 
try to delete a set of 
edges which 
will isolate either a single  vertex or a  
complete graph from $G$ in the hope of obtaining 
a simpler $(e)$-positive 
graph.  In 
order to see how this procedure will affect 
$Y_{G}$, we use the following lemma.

\ble                                             
  \label{secondary}
Given a graph, $G$ on $d$ vertices let
$H=G\uplus K_{m}$ where the vertices in $K_m$ are
labeled $v_{d+1},v_{d+2},\ldots,v_{d+m}$.  
If $Y_{G}=\ds{\sum_{\si \in 
\Pi_{d}}c_{\si}e_{\si}}$, then $Y_{H}= 
\ds{\sum_{\si \in 
\Pi_{d}}c_{\si}e_{\sigma/d+1,d+2, \ldots, 
d+m}}$.
\ele

{\bf Proof.} From the labeling of $H$ we have 
\begin{eqnarray*} Y_{H}&=&Y_{G}\* 
e_{[m]}\\
       &=&\sum_{\sigma \in 
\Pi_{d}}c_{\sigma}e_{\sigma}e_{[m]}\\
      &=& \sum_{\sigma \in 
\Pi_{d}}c_{\sigma}e_{\sigma/d+1,d+2, 
\ldots, 
d+m}.  \\ \end{eqnarray*}

 \vspace{-.2 in}\hfill $\Qed$

This result  suggests  we use 
the natural notation $G/v_{d+1}$ for the graph  
$G\biguplus \{v_{d+1}\}$.  We are now in a 
position to prove the $(e)-$positivity of paths.

\bpr  \label{paths} For all $d\geq 1$,
$Y_{P_{d}}$  is $(e)-$positive.
\epr

{\bf Proof.} We proceed by induction, having 
labeled  $P_{d}$ so that the edge set is 
$E(P_{d})=\{v_{1}v_{2}, v_{2}v_{3}, \ldots, 
v_{d-1}v_{d}\}. $
If $d=1$, then we have $Y_{P_{1}} = e_{1}$ and 
the proposition is clearly true.

So we assume by induction that 
$$Y_{P_{d}}\con_{d} \sum_{\pta \subseteq
\Pi_{d}}c_{\pta}e_{\pta},$$
 where $c_{\pta}\geq 0$ for all $\pta \in 
\Pi_{d}$.  
{}From the Deletion-Contraction Recurrence applied to $e=v_d v_{d+1}$,
%
%
Corollary \ref{dord+1} and  Lemma 
\ref{secondary}, we
see that
\begin{eqnarray*} \label{indp}
Y_{P_{d+1}}&=& Y_{P_{d}/v_{d+1}}- Y_{P_{d}} \ind 
\\
&\con_{d+1}&
 \sum_{\pta \subseteq 
\Pi_{d}}c_{(\tau)}e_{(\tau/d+1)}- 
\sum_{\pta \subseteq \Pi_{d}}c_{\pta}e_{\pta} 
\ind \\
&\con_{d+1}& 
\sum_{\pta \subseteq \Pi_{d}}c_{\pta}\l 
1-\frac{1}{|B_{\tau}|}\r e_{(\tau/d+1)}
+
\sum_{\pta \subseteq 
\Pi_{d}}{c_{\pta}\over{|B_{\tau}|}}e_{(\tau 
+(d+1))}.\\
\end{eqnarray*} 

Since we know that $c_{\pta}\geq 0$, and 
$|B_{\tau}| \geq 1$ for all $\tau$, this 
completes the induction step and the
 proof. \hfill $\Qed$

In the commutative context we will say that the 
symmetric function $X_{G}$ is $e$-{\em positive} 
if all the 
coefficients in 
the expansion of the elementary symmetric 
functions are non-negative.  
Clearly
$(e)$-positivity results for $Y_G$
specialize to $e$-positivity results
for $X_G$.  

\begin{corollary}
$X_{P_{d}}$ is $e-$positive. \hfill $\Qed$
\end{corollary}

One would expect the $(e)$-expansions for cycles and paths to be
related as is shown by the next proposition.
For labeling purposes, however, 
we first need a lemma which follows 
easily from the Relabeling Proposition.

\ble \label{fixed}  If $\gamma \in {\cal S}_{d}$ 
fixes $d$, then $Y_{\gamma(G)}\con_{d} Y_{G}. 
\hfill \Qed $ \ele


\bpr  For all $d\geq 1$, 
if $$Y_{P_{d}}\con_{d} \sum_{(\tau)}
c_{\pta}e_{\pta},\mbox{ then } 
Y_{C_{d+1}}\con_{d+1} 
\sum_{(\tau)} c_{\pta}e_{(\tau +(d+1))},$$ where we have 
labeled the graphs so $E(P_{d})=\{v_{1}v_{2}, 
v_{2}v_{3}, \ldots, 
v_{d-1}v_{d}\}$ and $E(C_{d+1})=\{v_{1}v_{2}, 
v_{2}v_{3}, \ldots, 
v_{d-1}v_{d},v_{d}v_{d+1},v_{d+1}v_{1}\}$.
\epr

{\bf Proof.}  We proceed by induction on $d$.
If $d=1$, then $Y_{P_{1}}=e_{[1]}$ and 
$Y_{C_{2}} = e_{[2]}$,  so the 
proposition holds for $d=1$.

For the induction step, we assume that 
$$Y_{P_{d-1}}\con_{d-1} \sum_{(\tau)} 
c_{\pta}e_{\pta}$$ and  also  that $$ 
Y_{C_{d}}\con_{d} \sum_{(\tau)} c_{\pta}e_{\ptap}.$$
  We notice that if $e=v_{d}v_{d+1}$, then 
$C_{d+1} -e$ does not have the standard labeling 
for paths.  But if we 
let 
$\ga=(d+1)(1,d)(2,d-1)\cdots (\left\lfloor 
\frac{d+1}{2} 
\right\rfloor , \left\lceil \frac{d+1}{2} 
\right\rceil)$ then we can
 use the Deletion-Contraction Recurrence to get 
$$Y_{C_{d+1}}= 
Y_{\gamma(P_{d+1})}-Y_{C_{d}}\ind.$$
However, since $d+1$ is a fixed point for 
$\gamma,$ Lemma \ref{fixed} allows us to deduce 
that
$$Y_{C_{d+1}}\con_{d+1}  
Y_{P_{d+1}}-Y_{C_{d}}\ind.$$

In the proof of Proposition \ref{paths} we saw 
that  $$Y_{P_{d+1}}=Y_{P_{d}/v_{d+1}} - 
Y_{P_{d}}\ind.$$
Combining these two equations gives 
%
%
%
\begin{equation} \label{a} 
Y_{C_{d+1}}\con_{d+1}Y_{P_{d}/v_{d+1}} - 
Y_{P_{d}}\ind -Y_{C_{d}}\ind.
\end{equation}
The demonstration of Proposition \ref{paths} 
also showed us that  
\begin{equation}  \label{b} Y_{P_{d}} \con_{d} 
\sum_{(\tau)} \l \l c_{\pta} 
-\frac{c_{\pta}}{|B_{\tau}|}\r e_{\ptas} 
+\frac{c_{\pta}}{|B_{\tau}|}
e_{\ptap}\r .\end{equation} 
Applying Corollary \ref{dord+1} and Lemma 
\ref{secondary} yields \begin{eqnarray*}
Y_{P_{d}}\ind  \con_{d+1}& \ds{\sum_{(\tau)}} & 
\left[ \l 
c_{(\tau)}-\frac{c_{(\tau)}}{|B_{\tau}|}\r 
e_{(\tau/d/d+1)} -\l 
c_{(\tau)}-\frac{c_{(\tau)}}{|B_{\tau}|}\r  
e_{(\tau/d,d+1)}\right. \\
&+& \left. 
\frac{c_{(\tau)}}{|B_{\tau}|(|B_{\tau}|+1)}  
e_{(\tau + (d)/d+1)} - 
\frac{c_{(\tau)}}{|B_{\tau}|(|B_{\tau}|+1)}  
e_{(\tau + (d)+ (d+1))}\right] 
\end{eqnarray*} 
and
$$Y_{P_{d}/v_{d+1}} \con_{d+1} \sum_{(\tau)} \l 
c_{(\tau)}-\frac{c_{(\tau)}}{|B_{\tau}|}\r 
e_{(\tau /d/d+1)} 
+\frac{c_{(\tau)}}{|B_{\tau}|}  
e_{(\tau+(d)/d+1)}.$$
By the induction hypothesis, 
\begin{eqnarray*}Y_{C_{d}}\ind &\con_{d+1}& 
\sum_{(\tau)} c_{\pta}e_{\ptap}\ind \\
 &\con_{d+1}& \sum_{(\tau)} 
\left(\frac{c_{\pta}}{|B_{\tau}|+1}
e_{\ptaps} - \frac{c_{\pta}}{|B_{\tau}| 
+1}e_{\ptapp}\right). \end{eqnarray*}

Plugging these  expressions for $ Y_{P_{d}/d+1}, 
 Y_{P_{d}}\ind,$ and $ Y_{C_{d}}\ind$  into 
equation (\ref{a}), 
grouping 
the terms according 
%
%
%
%
  to type, and simplifying gives
$$ Y_{C_{d+1}}\con_{d+1} \sum_{(\tau)} \left( 
c_{\pta}-\frac{c_{\pta}}{|B_{\tau}|}\right)
e_{\ptasp} +
\frac{c_{\pta}}{|B_{\tau}|}e_{\ptapp}.$$  This 
corresponds to the expression in equation 
(\ref{b}) for $Y_{P_{d}}$ 
in 
exactly the desired manner, and so we are done. 
$\hfill \Qed$ 

{From} the previous proposition and the fact that $Y_{C_1}=0$ we get an
immediate corollaries.
\bpr
For all $d\geq 1$, $Y_{C_{d}}$ is 
$(e)-$positive. \hfill $\Qed$
\epr
 
\bco 
For all $d\geq 1$, $X_{C_{d}}$ is $e-$positive. 
\hfill  $\Qed$
\eco

We are also able to use our recurrence to show 
the $(e)$-positivity of complete graphs with one 
edge removed.

\bpr
For  $d\geq 2$, if $e=v_{d-1}v_{d}$ then 
$$Y_{K_{d}-e} \con_{d} \frac{d-2}{d-1}e_{([d])} 
+\frac{1}{d-1}e_{([d-1]/d)}.$$  
\epr
%
{\bf Proof.}  Consider the complete graph 
$K_{d}$ and apply deletion-contraction to the 
edge $e=v_{d-1}v_{d}$.    
Together 
with Corollary \ref{dord+1} this will give us

\begin{eqnarray*}         e_{[d]}&=& Y_{K_{d}}\\
                                 &=& 
Y_{K_{d}-e}-Y_{K_{d-1}}\ind\\
&=&Y_{K_{d}-e}-e_{[d-1]}\ind \\
                                 &\con_{d}& 
Y_{K_{d}-e}-\frac{1}{d-1}e_{([d-1]/d)} + 
\frac{1}{d-1}e_{([d])}.\\
                                 \end{eqnarray*}
                                 Simplifying 
gives the result.  \hfill $\Qed$

This also immediately specializes.

\bco
For  $d \geq 2$, 
$$X_{K_{d}-e}=d(d-2)(d-2)!e_{d}+(d-2)!e_{(d-1,1)
}.$$
\vspace{-.2 in} \hfill $\Qed$
\eco

\section{The ({\bf 3+1})-free  Conjecture}
One of our original goals in setting up this 
inductive machinery was to make progress on 
 the ({\bf 3+1})-free  Conjecture of 
Stanley and Stembridge, which we now state.
Let {\bf a+b}
be the poset which is a disjoint union of an 
$a$-element chain and a $b$-element chain. The 
poset $P$ is said to be 
({\bf 
a+b})-{\em free} if it contains no induced 
subposet isomorphic to {\bf a+b}.  Let $G(P)$ 
denote the incomparability graph of $P$ 
whose vertices are the elements of $P$ with
an edge $uv$ whenever $u$
and $v$ are incomparable in $P$.  The ({\bf 
3+1})-free Conjecture  of Stanley and Stembridge 
\cite{stembridge} states:  
\begin{conjecture} If $P$ is ({\bf 
3+1})-free, 
then $X_{G(P)}$ is $e$-positive. 
\end{conjecture}
Gasharov~\cite{gasharov} has demonstrated the weaker result that
$X_{G(P)}$ is $s$-positive, where $s$ refers to the Schur functions.

A subset of the ({\bf 3+1})-free graphs is the 
class of {\em indifference graphs}.   They are 
characterized 
\cite{stanley} 
as having vertices and edges 
$$
\mbox{  $V=\{v_1,\ldots,v_d\}$
and   $E=\{v_iv_j: i,j \mbox{ belong to some }[k,l]\in 
\cal{C}\}$,}
$$ 
where $\cal{C}$ is  a collection of intervals 
$[k,l]=\{k,k+1,\ldots, l\} \subseteq  [d]$.  We 
note that without loss of 
generality, we can assume no interval in the 
collection is properly contained in any other.  
These correspond to incomparability graphs 
of posets which 
 are both 
({\bf 3+1})-free and ({\bf 2+2})-free.

Indifference graphs have a nice inductive structure that should make
it possible to apply our deletion-contraction techniques.  Although we
have not been able to do this for the full family,  we are able to
resolve a special case.   For any composition of $n$,
%
%
  $\al = (\al_{1},\al_{2}, \ldots, 
\al_{k}), $ let 
$\tilde{\al_{i}}=\sum_{j\leq i}\al_{j}$.  
A 
{\it $K_{\al}$-chain} is the  indifference graph 
using the collection of intervals 
$\{[1,\tilde{\al_{1}}],[\tilde{\al_{1}},
\tilde{\al_{2}}],\ldots, 
[\tilde{\al_{k-1}},\tilde{\al_{k}}]\}.$  
This is  just a string of complete graphs, whose 
sizes are given by 
the 
parts of $\al$, which are attached to one 
another sequentially at single vertices.  We 
notice that the 
$K_{\al}$-chain for $\al = 
(\al_{1},\al_{2}, \ldots, \al_{k})$ 
can be obtained from the 
$K_{\tau}$-chain for $\tau = 
(\al_{1},\al_{2}, \ldots, 
\al_{k-1})$ by attaching the graph 
$K_{\al_{k}}$ to  
its last vertex.

We will be able to handle this  type of 
attachment  for any graph $G$ with vertices $\{ 
v_{1},v_{2},\ldots,v_{d}\}.$  
Hence, we define $G+K_{m}$ to be the graph with 
$$\mbox{$V(G+K_{m})=V(G)\cup \{ 
v_{d+1},\ldots,v_{d+m-1} \}$}$$ and 
$$\mbox{ $E(G+K_{m})=E(G)\cup \{ e=v_{i}v_{j}: 
i,j \in [d,d+m-1] \}.$}$$   Using 
deletion-contraction techniques, we 
are 
able to exhibit the relationship between the 
$(e)$-expansion of $G+K_{m}$ and the  
$(e)$-expansion of 
$G$.  However,  we will also need some more 
notation.  For $\pi \in \Pi_{d}$, we  let
 $\pi +i $ denote the partition given by $\pi$ 
with the additional $i$ 
elements $d+1,d+2, \ldots, d+i$  added to 
$B_{\pi}$.  This is in contrast to  $ \pi +(i)$, 
which denotes the 
partition given by $\pi$ with the element $i$ 
inserted into 
$B_{\pi}$. We denote the falling factorial by
$$\ipr{m}_{i}\defn m(m-1)\cdots (m-i+1)$$ and 
the rising factorial by
$$(m)_{i} \defn m(m+1)\cdots (m+i-1).$$

We begin studying the behavior of 
$Y_{G+K_{m}}\ind_{d}^{d+j}$ with two lemmas.

\ble  If $1 \leq j < k \leq m$, then 
$Y_{G+K_{m}}\ind_{d}^{d+j}\con_{d}
Y_{G+K_{m}}\ind_{d}^{d+k}.$   \label{jk} \ele

{\bf Proof. } For any partition, $\pi$ of 
$[d+m-1]$,, define $\delta_{d+j}(\pi)$ to be the 
partition  obtained by inserting 
the element $d+j$ into the block of $\pi$ 
containing $d$, and adding one to each element 
of $\pi$ which 
is at least $d+j$.  It is easy to see 
that 
$m_{\pi}\ind_{d}^{d+j}=m_{\delta_{d+j}(\pi)}.$
Similarly, 
$m_{\pi}\ind_{d}^{d+k}=m_{\delta_{d+k}(\pi)}.$  
Now consider the permutation $\gamma = 
(d+j,d+k,d+k-1,d+k-2,\ldots,d+j+1)$ constructed so 
that $\gamma \comp 
(m_{\pi}\ind_{d}^{d+j})=m_{\pi}\ind_{d}^{d+k}.$  
This 
implies that  $\gamma \comp (Y_{G+K_{m}} 
\ind_{d}^{d+j})=Y_{G+K_{m}}\ind_{d}^{d+k}\hspace
{-1mm}.$ 
Noticing that $d$ is a fixed point of $\gamma$ 
(so that Lemma \ref{fixed} applies) will complete 
the proof. \hfill $\Qed$

\ble\label{prelim}

If $G$ is a graph on $d$ vertices with $$Y_{G} 
\con_{d}
 \sum_{(\pi)\subseteq 
\Pi_{d}}c_{(\pi)}e_{(\pi)},$$ then $$
  Y_{G+K_{m}}\hspace{-1mm}\ind_{d}^{d+m} 
\con_{d+m}\hspace{-1mm}\sum_{(\pi)}
\sum_{i=0}^{m-1} 
\frac{c_{(\pi)}\ipr{m-1}_{i}\left[ 
 e_{(\pi+i/d+i+1,\ldots,d+m)}-e_{( \pi + i 
+(d+m)/d+i+1,\ldots, d+m-1)}\right ] 
}{(b)_{i+1}}, $$
%
%
%
 \n  where $b=|B_{\pi}|$.

\ele

{\bf Proof.}  We prove the lemma by induction on 
$m$.  The case $m=1$ 
is merely a restatement of Corollary 
\ref{dord+1}.  So we may 
 assume this lemma is true for 
$Y_{G+K_{m}}\ind_{d}^{d+m}$, and 
proceed to prove it for
$Y_{G+K_{m+1}}\ind_{d}^{d+m+1}.$

{From} Lemma \ref{jk}, it follows that  for   
$1\leq j \leq m$, 
we have $$Y_{G+K_{m}}\ind_{d}^{d+j}\h 
\ind_{d}^{d+m+1} \con_{d+m+1} 
Y_{G+K_{m}}\ind_{d}^{d+m}\h \ind_{d}^{d+m+1}.$$

Now, from $G+K_{m+1}$ we may  delete  the edge
set  $\ds{\{v_{d}v_{d+j}: 1\leq j \leq m\}}$   
and combine all the terms $Y_{G}\ind_{d}^{d+j}\h 
\ind_{d}^{d+m+1}$ for $1\leq j \leq m$ to obtain 
\begin{eqnarray*}Y_{G+K_{m+1}}\ind_{d}^{d+m+1} 
&\con_{d+m+1}& Y_{G\uplus K_{m}}\ind_{d}^{d+m+1} 
-mY_{G+K_{m}}\ind_{d}^{d+m}\h \ind_{d}^{d+m+1}\\
 & \con_{d+m+1}& Y_{G\uplus 
K_{m}}\ind_{d}^{d+m+1} 
-mY_{G+K_{m}}\ind_{d}^{d+m}\h 
\ind_{d+m}^{d+m+1}.
 \end{eqnarray*}

{From} this point on, we need only concern 
ourselves with the clerical 
details, making  sure that everything matches up 
properly.  We can see from Lemma \ref{main}, 
Lemma \ref{secondary} 
and the 
original hypothesis on $Y_{G}$,  that 
\begin{equation}                 \label{mess1}
Y_{G \uplus K_{m}}\ind_{d}^{d+m+1} \con_{d+m+1} 
\sum_{(\pi)} \frac{c_{(\pi)}}{b}\l 
e_{(\pi_{1})}-
e_{(\pi_{2})}\r .
\end{equation}
where

$\begin{array}{l}
 \pi_{1}=\pi/d+1,\ldots,d+m/d+m+1,   \\
\pi_{2}=\pi+(d+m+1)/d+1,\ldots,d+m .\\
\end{array}$


\vspace{.2 in}

    Similarly, the induction hypothesis shows 
\begin{equation}                       
\label{mess2}
mY_{G+K_{m}}\ind_{d}^{d+m}\h \ind_{d+m}^{d+m+1} 
\con_{d+m+1} \newline 
\sum_{(\pi)}  \sum_{i=0}^{m-1}\frac{c_{(\pi)}m 
\ipr{m-1}_{i}}{(b)_{i+1}}
\l    \frac{e_{(\pi_{3})} - e_{(\pi_{4})}}{m-i} 
- \frac{ e_{(\pi_{5})}-e_{(\pi_{6})}}{b+i+1}\r
\end{equation}

\n where

$\begin{array}{l}
\pi_{3}=\pi+i/d+i+1,\ldots,d+m/d+m+1,    \\
\pi_{4}=\pi+i/d+i+1,\ldots,d+m+1,\\
\pi_{5}=\pi+i+(d+m)/d+i+1,\ldots,d+m-1/d+m+1, \\
\pi_{6}=\pi+i+(d+m)+(d+m+1)/d+i+1,\ldots,d+m-1. 
\\
\end{array}$

\vspace{.2 in}

Simplifying the  terms  and combining both 
equations (\ref{mess1}) and 
(\ref{mess2}) gives   
$$\begin{array}{l}
\ds{ Y_{G+K_{m+1}}\ind_{d}^{d+m+1} \con_{d+m+1}} 
\\
\\
\qquad \ds{\sum_{(\pi)} c_{(\pi)}\l \frac{ 
e_{(\pi_{1})}-e_{(\pi_{2})}}{b} 
   - \sum_{i=0}^{m-1}\frac{\l e_{(\pi_{3})} 
-e_{(\pi_{4})}\r\ipr{m}_{i}}{(b)_{i+1}}  
+\sum_{i=0}^{m-1}\frac{\l 
e_{(\pi_{5})} -e_{(\pi_{6})}\r 
\ipr{m}_{i+1}}{(b)_{i+2}} \r.} \\
 \end{array}$$

\n Note  that modulo $d+m+1$ we have 

\vspace{.1in}

$\begin{array}{l}
(\pi_{5})=(\pi+i+1/d+i+2,\ldots,d+m/d+m+1)  
\mbox{ and}\\
(\pi_{6})=(\pi +i+1+(d+m+1)/d+i+2,\ldots,d+m).
  \end{array}$

\vspace{.2 in}

\n So by shifting indices and simplifying, we 
obtain  
 $$\begin{array}{l}
 \ds{ Y_{G+K_{m+1}}\ind_{d}^{d+m+1} 
\con_{d+m+1}}\\
 \\ 
\qquad \ds{\sum_{(\pi)}  
\sum_{i=0}^{m}\frac{c_{(\pi)}\ipr{m}_{i} \left[ 
e_{(\pi+i/d+i+1,\ldots,d+m+1)} 
-e_{(\pi+i+(d+m+1)/d+i+1,\ldots,d+m)}
 \right]   }{(b)_{i+1}}},\end{array} 
$$
which completes the induction step and the 
proof. \hfill $\Qed$


This lemma is useful because it will helps us to find 
an explicit formula for 
$Y_{G+K_{m+1}}$ in terms of $Y_{G}$.  Once this formula is in hand, 
it will be easy to verify 
that if $G$ is $(e)$-positive, then so is 
$G+K_{m+1}$.  To complete the induction step in 
establishing this  formula, 
we will need the following observation.

\ble \label{combine}
For any graph $G$ on $d$ vertices, and $(i,j)\in 
{\cal S}_{d}$,  
$$Y_{G}\ind_{i}^{d+1}\con_{d+1}Y_{(i,j)(G)}\ind_
{j}^{d+1}.$$
\ele

{\bf Proof: } For any $\pi \in \Pi_{d}$, we see 
directly from the definitions for induction and 
the symmetric group 
action that 
$$(i,j)\comp(m_{\pi}\ind_{i}^{d+1})=m_{(i,j)\pi}
\ind_{j}^{d+1}.$$
It follows that $$(i,j)\comp 
(Y_{G}\ind_{i}^{d+1})=Y_{(i,j)(G)}\ind_{j}^{d+1}
.$$  Since $d+1$ is a fixed 
point of $(i,j)$, Lemma \ref{fixed} gives the 
result. \hfill $\Qed$

We now give the formula for $Y_{G+K_{m+1}}$ 
in terms of $Y_{G}$.

\ble  \label{add}
If  $m\geq 1$, and $$  Y_{G} \con_{d} 
\sum_{(\pi)\subseteq 
\Pi_{d}}c_{(\pi)}e_{(\pi)},$$ then $$ 
Y_{G+K_{m+1}}\con_{d+m} 
 \sum_{(\pi)\sbe \Pi_{d}}  \sum_{i=0}^{m-1} 
\frac{c_{(\pi)}\ipr{m-1}_{i} 
}{(b)_{i+1}} \left[  (b-m+i)
 e_{(\hat{\pi})} +(i+1)e_{(\ol{\pi})} \right]  
$$ 
where $b=|B_{\pi}|$ and

$\begin{array}{l}
 \hat{\pi}= \pi +i/d+i+1,\ldots , d+m,\\
 \ol{\pi}=\pi + i +(d+m)/ d+i+1,\ldots,d+m-1.\\
\end{array}$
\ele

{\bf Proof.}  We induct on $m$.  If $m=1$, then 
$Y_{G+K_{2}} 
=Y_{G\uplus K_{1}} - Y_{G}\ind_{d}^{d+1}.$
This shows that $$Y_{G+K_{2}} \con_{d+1} 
\sum_{(\pi)} 
\l \frac{c_{(\pi)}(b-1)}{b}e_{(\pi/d+1)} + 
\frac{c_{(\pi)}}{b}e_{(\pi+(d+1))}\r,$$ which 
verifies the base case.

To begin the induction step, we repeatedly 
utilize the 
Deletion-Contraction Recurrence to delete the 
edges $v_{d+i}v_{d+m+1}$ for $0 \leq i \leq m$, 
and obtain 
\begin{equation} \label{inductionstep} 
Y_{G+K_{m+2}}\con_{d+m+1} Y_{G+K_{m+1}\uplus 
v_{d+m+1}}-mY_{G+K_{m+1}}\ind_{d+m}^{d+m+1} -  
Y_{G+K_{m+1}}\ind_{d}^{d+m+1}. \end{equation}
Note that we are able to combine all the terms 
from  $Y_{G+K_{m+1}}\ind_{d+i}^{d+m+1}$ for $1\le i\le m$ using
Lemma \ref{combine}, 
since in these cases the
necessary permutation exists.

We now expand
 each of the terms in equation 
(\ref{inductionstep}). For the first, using 
Lemma \ref{secondary}, $$Y_{G+K_{m+1}\uplus 
v_{d+m+1}}\con_{d+m+1}
\sum_{(\pi)}  \sum_{i=0}^{m-1} 
\frac{c_{(\pi)}\ipr{m-1}_{i}
 }{(b)_{i+1}} \left[ (b-m+i) e_{(\pi_{1} )} + 
(i+1)e_{(\pi_{2} )}   \right], $$
where
$$\begin{array}{l}
  \pi_{1}=\pi +i/d+i+1,\ldots , d+m/d+m+1 , \\
  \pi_{2}=\pi + i +(d+m)/ d+i+1,\ldots,d+m-1/d+m+1.
\end{array}$$

\vspace{.2 in}

\n For the second term, using Corollary 
\ref{dord+1}, we have
$$\begin{array}{l} 
\ds{mY_{G+K_{m+1}}\ind_{d+m}^{d+m+1} 
\con_{d+m+1}} \\ 
\\ 
 \qquad 
\ds{\sum_{(\pi)}\sum_{i=0}^{m-1}\frac{c_{(\pi)}
\ipr{m}_{i+1}}{(b)_{i+1}}
\left[ \frac{b-m+i}{m-i} \l e_{(\pi_{1 
})}-e_{(\pi_{3})}\r + \frac{i+1}{b+i+1}\l 
e_{(\pi_{2})}-e_{(\pi_{4})} \r  
\right],}\\ \end{array}$$
where
$$\begin{array}{l}
         \pi_{3}=\pi+i/d+i+1,\ldots,d+m+1, \\
\pi_{4}=\pi+i+(d+m)+(d+m+1)/d+i+1,\ldots,d+m-1.\\                  
\end{array}
$$

\vspace{.2 in}

\n And finally,  using Lemma \ref{prelim}, 

$$Y_{G+K_{m+1}}\ind_{d}^{d+m+1}\con_{d+m+1}\\
\sum_{(\pi)}\sum_{i=0}^{m} 
\frac{c_{(\pi)}\ipr{m}_{i}}{(b)_{i+1}} 
\l e_{(\pi_{3})} - e_{(\pi_{5})}\r $$

\n where 

$\begin{array}{l}
        \pi_{5}=\pi+i+(d+m+1)/d+i+1,\ldots,d+m 
.\\
\end{array}$

\vspace{.2 in}

\n Grouping the terms appropriately and shifting 
indices where needed gives

$$Y_{G+K_{m+2}}\con_{d+m+1} 
\sum_{(\pi)}\sum_{i=0}^{m}\frac{c_{(\pi)}
\ipr{m}_{i} \left[ 
(b-(m+1)+i)
e_{(\pi_{3})} + 
(i+1)e_{(\pi_{5})}\right]}{(b)_{i+1}} .$$

This completes the induction step and the proof. 
\hfill $\Qed$

Examining this lemma, we can see that in ${\ds 
Y_{G+K_{m+1}}}$ 
we have the same sign on all the coefficients as 
we had in ${\ds Y_{G}}$, 
with the possible exception of the terms where 
$b<m-i$.  But it is easy to see that in this 
case we  have 
$$e_{(\pi+i/d+i+1,\ldots,d+m)}\con_{d+m}
e_{(\pi+m-i-b-1+(d+m)/d+m-i-b,\ldots,d+m-1)}.$$
This means that in the expression for ${\ds 
Y_{G+K_{m+1}}}$ as a sum 
over  congruence classes modulo $d+m$, we can 
combine the coefficients on
 these terms.  And so upon simplification, the 
coefficient
 on $e_{(\pi+i/d+i+1,\ldots,d+m)}$ will be:
 
%

$$\l 
\frac{(b-m+i)\ipr{m-1}_{i}}{(b)_{i+1}}+\frac{(m-
i-b)\ipr{m-1}_{m-i-b-1}}
{(b)_{m-i-b}}\r c_{(\pi)},$$
 where $c_{(\pi)}$ is the coefficient on 
$e_{(\pi)}$ in $Y_{G}$.

Adding these fractions by finding a common 
denominator, we see that this is  
actually zero, which gives us the next result.

\bth \label{chain}
If $Y_{G}$ is  $(e)$-positive, then 
$Y_{G+K_{m}}$ is also $(e)$-positive.
\hfill $\Qed$
\eth

Notice that Proposition \ref{paths} follows 
easily from Theorem \ref{chain} and induction, 
since for paths 
$P_{m+1}=P_{m}+K_{2}$.  As a more general result 
we have the following corollary.

\bco
If $G$ is a $K_{\al}$-chain, then $Y_{G}$ is 
$(e)$-positive.  Hence, $X_{G}$ is also 
$e$-positive. \hfill $\Qed$
\eco

We can also describe another  class of 
$(e)$-positive graphs.
We define a {\em diamond} to be the indifference 
graph on the collection of  intervals  
$\{[1,3],[2,4]\}.$  So a 
diamond 
consists of two $K_{3}$'s sharing a common edge. 
 Then the following holds.


\bth
Let $D$ be a  diamond.  If $G$ is 
$(e)$-positive, then so is $G+D$.
\eth

{\bf Proof.}   The proof of
 this result is analogous to the proof for the 
case of 
$G+K_{m}$, and so is omitted. \hfill
$\Qed$

 \section{Comments and open questions}
 
We will end with some questions raised by this work.  
We hope they will stimulate future research.

\medskip

(a)  Obviously it would be desirable to find a way to use
deletion-contraction to prove that indifference graphs are
$e$-positive (or even demonstrate the full ({\bf 3+1})-Free
Conjecture).  The reason that it becomes difficult to deal with the case
where the last two complete graphs overlap in more than one vertex
is because one has to keep track of all ways the intersection could be
distributed over the block sizes of an $e_\pi$.  Not only is the
bookkeeping complicated, but it becomes harder to find groups of
coefficients that will sum to zero.
 
Another possible approach is to note that if $G$ is an
indifference graph, then for the edge $e=v_{k}v_{d}$ (where $[k,d]$ is
the  interval) both $G\setm e$ and $G/e$ are indifference 
graphs.  Furthermore $G\setm e$ is obtained from $G/e$ by attaching a
$K_{d-k}$ so that it intersects in all but one vertex with the final
$K_{d-k}$ of $G\setm e$.  Unfortunately, the 
relationship between the 
coefficients in the $(e)$-expansion of $Y_{G\setm e}$ 
and $Y_{G/e}\ind$ does not seem to be very simple. 

\medskip 
 
(b)  Notice that if $T$ is a tree on 
$d$ vertices, we have ${\cal 
X}_{T}(n)=n(n-1)^{d-1}.$  Since $X_{G}$ is 
a generalization of the chromatic polynomial, it 
might be reasonable to suppose that it also is 
constant on trees with 
$d$ vertices.  This is far from the case!  In 
fact, it has been verified up to $d=9$ 
\cite{tim} that, for non-isomorphic trees $T_1,T_2$ 
we have $X_{T_1}\neq X_{T_2}$.
This leads to the following question posed by 
Stanley.
\begin{question}[\cite{stan}]
Does $X_{T}$ distinguish between non-isomorphic trees?
 \end{question}

We should note that  the answer to this question is definitely ``yes''
for $Y_T$.  In fact more is true.
\bpr
The function $Y_G$ distinguishes between all graphs $G$
with no loops or multiple edges.
\epr
 
{\bf Proof.}  We know from Proposition \ref{exp}
that $Y_{G}=\sum_{P}m_{\pi(P)}$ for the stable 
partitions $P$.  
Construct the graph $H$ with vertex set 
$V(G)=\{v_{1},v_{2},\ldots,v_{d}\}$ and edge set 
 $E(H)=\{v_{i}v_{j}| \mbox{ 
there exists a $\pi(P)$ such that $i,j$ are in 
the same block of $\pi(P)$}\}.$  Since $\pi(P)$ 
comes from a stable 
partition $P$ of $G$, $v_{i}$ and $v_{j}$ are in 
the same block of some $\pi(P)$ if and only if 
there is no edge 
$v_{i}v_{j}$ in $G$.  Hence the graph $H$ 
constructed is the (edge) complement of $G$ and 
so we can recover $G$ from 
$H$. \hfill $\Qed$
 
Of course we can have $Y_G\neq Y_H$ but $X_G=X_H$.
So a first step towards answering Stanley's question might be to see
if $Y_T$ still distinguishes trees under congruence.
It seems reasonable to expect to investigate this using our 
deletion-contraction techniques since trees are 
reconstructible from their leaf-deleted 
subgraphs  \cite{ed}.  We 
proceed in the following manner.
 
 If $T_{1}\not\cong T_{2}$ then by the 
reconstructibility of trees there must exist 
labelings of these trees so that 
$v_{d}$ is a leaf of $T_{1}$, $\tilde{v}_{d}$ is 
a leaf of $T_{2}$ and $T_{1}-v_{d} \not\cong 
T_{2}-\tilde{v}_{d}$.  By 
induction we will have 
$Y_{T_{1}-v_{d}}\not\con_{d-1}Y_{T_{2}-\tilde{v}_{d}}$, and
consequently,  $Y_{T_{1}-v_{d}}\not\con_{d}Y_{T_{2}-\tilde{v}_{d}}.$
Furthermore,  our recurrence gives  
$$\begin{array}{l}
Y_{T_{1}}=Y_{T_{1}-v_{d}/v_{d}}-Y_{T_{1}-v_{d}}
\ind \\
Y_{T_{2}}=Y_{T_{2}-\tilde{v}_{d}/\tilde{v}_{d}}-
Y_{T_{2}-\tilde{v}_{d}}\ind. \\ 
\end{array}$$
One now needs to investigate what sort of cancelation occurs to see if
these two differences could be equal or not.  Concentrating on a term
of a particular type could well be the key.

\medskip

(c)  It would be very interesting to develop a wider theory of
symmetric functions in noncommuting variables.  The only relevant
paper of which we are aware is Doubilet's~\cite{doubilet} where he
talks more generally about functions indexed by set partitions, but
not the noncommutative case per se.  His work is dedicated to finding
the change of basis formulae between 5 bases (the three we have
mentioned, the complete homogeneous basis, and the so-called forgotten
basis which he introduced).  However, there does not as yet seem to be
any connection to representation theory.  In particular, there is no
known analog of the Schur functions in this setting.


\begin{thebibliography}{15}

\bibitem{bs} A. Blass and B. Sagan, Bijective 
proofs of two broken circuit theorems, {\em J. 
Graph Theory} {\bf 10} 
(1986), 15--21.
\bibitem{tim} T. Chow, personal communication. 
\bibitem{doubilet} P. Doubilet, On the 
Foundations of Combinatorial Theory. VII: 
Symmetric Functions through 
the Theory of Distribution and Occupancy, {\em 
Studies in Applied Math.} {\bf 51} (1972), 
377--396. 
\bibitem{gasharov} V. Gasharov, Incomparability 
graphs of ({\bf 3}+{\bf 1})-free posets are 
$s$-positive, {\em 
Discrete Math.} {\bf 157} (1996), 193--197
\bibitem{thibon}I.\ M.\ Gelfand, D.\ Krob, A.\ 
Lascoux, B.\ Leclerc, V.\ Retakh, J.-I.\ Thibon, 
Noncommutative 
symmetric functions, preprint, 1994.

\bibitem{orient} D. Gebhard and B. Sagan, Sinks in acyclic
orientations of graphs, preprint, 1999.

\bibitem{gandz}  C.\ Greene and T.\ Zaslavsky, 
On the interpretation of Whitney numbers through 
arrangements of 
hyperplanes, zonotopes, non-radon partitions, 
and orientations of graphs, {\em Trans. Amer. 
Math. Soc.} {\bf 280} 
(1983), 97--126.
\bibitem{me}  D.\ Gebhard, Noncommutative 
symmetric functions and the chromatic 
polynomial, \emph{Conference 
Proceedings at the 7th International Conference 
on Formal Power Series and Algebraic 
Combinatorics}, 
Universit\'e de Marne-la-Vall\'ee 1995.
\bibitem{ed} F.\ Harary and E.\ Palmer, The 
reconstruction of a tree from its maximal 
subtrees, {\em Canad. J. Math.}, 
{\bf 18} (1966), 803--810.




\bibitem{??} R.\ P.\ Stanley, Acyclic 
orientations of graphs, {\em Discrete Math.} 
{\bf 5} (1973), 171--178.


\bibitem{nw:wgp} S. D. Noble and D. J. A. Welsh, A weighted graph
polynomial from chromatic invariants of knots, Annales l'Institut
Fourier 1999, to appear.


\bibitem{stan} R.\ P.\ Stanley, A symmetric 
function generalization of the chromatic 
polynomial of a graph, 
{\em Advances in Math.} {\bf 111} (1995), 
166--194.

\bibitem{stanley} R.\ P.\ Stanley, Graph 
colorings and related symmetric functions: ideas 
and applications, 
preprint.

\bibitem{stembridge} R.\ P.\ Stanley and J.\ 
Stembridge, On immanants of Jacobi-Trudi 
matrices and permutations with 
restricted position, {\em J. Combin. Theory (A)} 
{\bf 62} (1993), 261--279.
\bibitem{whitney} H.\ Whitney, A logical 
expansion in mathematics, {\em Bull. Amer. Mth. 
Soc.} {\bf 38} 
(1932), 572--579.
 


\end{thebibliography}
\end{document}